%% file: main.tex
\documentclass[journal]{IEEEtran} 

\usepackage{custom_cmd3}
\usepackage{longtable}
\usepackage{booktabs,subcaption,amsfonts,dcolumn}

\begin{document}

\title{Polyhedral study of a temporal rural postman problem: application in inspection of railway track without disturbing train schedules} 
\author{Somnath Buriuly, Leena Vachhani, Sivapragasam Ravitharan, Arpita Sinha, Sunita Chauhan
\thanks{S. Buriuly is a postdoctoral fellow with CoE-OGE, IIT Bombay, India}
\thanks{L. Vachhani and A. Sinha are with Systems and Control Engineering,
IIT Bombay, India}
\thanks{S. Ravitharan is with the Department of Mechanical and Aerospace Engineering, Monash University, Australia}
\thanks{S. Chauhan is with the Center for Equitable $\&$ Personalized Health, Plaksha University, India}}






\markboth{Manuscript for submission}%
{Buriuly \MakeLowercase{\textit{et al.}}: Polyhedral study of a temporal rural postman problem: application in inspection of railway track without disturbing train schedules} 

\maketitle

\begin{abstract}

\cora The Rural Postman Problem with Temporal Unavailability (RPP-TU) is a variant of the Rural Postman Problem (RPP) specified for multi-agent planning over directed graphs with temporal constraints. These temporal constraints represent the unavailable time intervals for each arc during which agents cannot traverse the arc. Such arc unavailability scenarios occur in non-disruptive routing and scheduling of the instrumented wagons for inspecting railway tracks, without disturbing the train schedules, i.e. the scheduled trains prohibit access to the tracks in the signal blocks (sections of railway track separated by signals) for some finite interval of time. \corb 

A three-index formulation for the RPP-TU is adopted from the literature. The three-index formulation has binary variables for describing the route information of the agents, and continuous non-negative variables to describe the schedules at pre-defined locations. 
A relaxation of the three-index formulation for RPP-TRU, referred to as Cascaded Graph Formulation (CGF), is investigated in this work. The CGF
has attributes that simplify the polyhedral study of time-dependent arc routing problems like RPP-TRU.
A novel branch-and-cut algorithm is proposed to solve the RPP-TU, where branching is performed over the service arcs. 
A family of facet-defining inequalities, derived from the polyhedral study, is used as cutting planes in the proposed branch-and-cut algorithm to reduce the computation time by up to $48\%$. 
Finally, an application of this work is showcased using a simulation case study of a railway inspection scheduling problem based on Kurla-Vashi-Thane suburban network in Mumbai, India. An improvement of $93\%$ is observed when compared to a Benders' decomposition based MILP solver from the literature. 
\end{abstract} 


\IEEEpeerreviewmaketitle

\section{Introduction}
\label{bncrpptu:sec:intro}

\subfile{sec-intro.tex}

\section{Preliminaries and Terminologies}
\label{bncrpptu:sec:defs}

\subfile{sec-prelims.tex}

\section{The Rural Postman Problem with Temporal Unavailability} 
\label{bncrpptu:sec:rpptucgf}

\subfile{sec-rpptu.tex}

\section{Proposed algorithm \cors for solving RPP-TU\corb }
\label{bncrpptu:sec:poly_algo}

\subfile{sec-algo.tex}


\section{Results and Discussion}
\label{bncrpptu:sec:results}

\subfile{sec-res.tex}

\section{Conclusion}
\label{bncrpptu:sec:conc}

We discuss our formulation for the Rural Postman Problem with Temporal Unavailability (RPP-TU), by modifying the replicated graph concept from existing literature. The construction of a replicated graph is conceptually an extension to the Arc Path Arc Sequence (APAS) formulation for multi-agent problems. 
In this work, we present a polyhedral study of the framework and determine families of valid facet-defining inequalities. A family of valid inequalities is implemented in a branch-and-cut algorithm, such that these inequalities act as cutting planes and improve the convergence time of our algorithm. Theorem \ref{thm:conn2b} doesn't require the set $S$ to be connected and it is not necessary for the set $\mathcal{A}_q$ to be entirely contained in $A(S) \cup \delta(S)$, thus making the separation algorithm simpler. 
In addition, we also propose a branching strategy over service arcs for the multi-agent routing problem. 

Extensive comparison of RPP-TU with the Time-dependent Rural Postman Problem (TDRPP) shows similarity in computation time for networks of the same size and similar attributes. This serves as a benchmark for our proposed methodology of formulating RPP-TU, and applying a branch-and-cut algorithm to solve for an exact solution. The comparison results also show up to $48\%$ improvement compared to the branch-and-cut algorithm without the proposed cutting planes. Since the RPP-TU is an unavailability-based routing and scheduling problem, we demonstrate an application using a simulation case study on Mumbai (India) based suburban railway network. The objective of this case study is to optimally inspect a specified section of track, where the tracks are available based on the daily suburban passenger train schedules. The resulting solution is an optimal plan for the inspection of specified railway track sections that don't disturb the existing train schedules by navigating only when the respective track sections are available. Observe that, in comparison with a Benders' decomposition based MILP solver from the literature, the proposed branch-and-cut algorithm shows $93\%$ reduction in computation time. 


\section*{Acknowledgment}
%
\corh This work is supported by IITB-Monash Research Academy. \corb 

\section*{Compliance with Ethical Standards} 

\noindent \textbf{Conflict of interest:} The authors declare that they have no conflict of interest. 

\noindent \textbf{Ethical approval:} This work does not contain any studies with human participants or animals performed by any of
the authors. 

\noindent \textbf{Informed consent:} Informed consent was obtained from all individual participants included in the study.

\printbibliography[title=References]

\section{Appendix}
\label{bncrpptu:sec:apdx}

\subfile{sec-apndx}




\end{document}

%% file: sec-intro.tex

\cora \IEEEPARstart{R}{outing} and scheduling problems arising in transportation networks often fall in the category of combinatorial optimization, hence \cors they are \cora modeled as discrete optimization problems. \corb 
Travelling Salesman Problem (TSP), Chinese Postman Problem (CPP) and Rural Postman Problem (RPP) are  popular combinatorial optimization problems defined for planning the routes of vehicles/agents in a transportation network. These combinatorial optimization problems are described using undirected graph $G := (V, E)$ and involve decisions to represent tours\footnote{a sequence of edges (describing movement in a network) that start and end at the same vertex (particular location in a network)} that fulfill various objectives; for example: in TSP, all the vertices $V$ must be visited only once (the graph is fully connected); in CPP, all the edges $E$ must be traversed at least once; and in RPP, some of the edges $E_R \subseteq E$ must be traversed at least once (a generalized version of CPP). 
We investigate a multi-agent variant of RPP on directed graphs with time-dependent availability of arcs (directed edges), called the Rural Postman Problem with Temporal Unavailability (RPP-TU). The proposed framework is suitable for modeling the routing and scheduling problems arising in railway networks because the track sections are separated by signals (called signal blocks) that will harbor only one locomotive at a given time. RPP-TRU is motivated by the railway track inspection planning problem using instrumented wagons so that the train schedules are not disrupted. This requires navigating through the available time intervals of the signal blocks. Conventionally, the inspections are conducted at night when the train services stop, see \cite{chen2023}. Hence, this work extends to the application where the inspection schedules cooperatively overlap with regular train schedules.

In combinatorial optimization problems like TSP, CPP, and RPP, the key decision is to determine the number of times an edge must be traversed. The order of traversal does not affect the optimality of the solution, and hence an optimal tour is easy to determine given the solution graph (a graph with the edges duplicated based on the optimal integer decision/solution, that represents the spatial solution without the traversal order). Once a solution graph is determined, a simple polynomial-time algorithm, to determine an Eulerian cycle, is sufficient to construct an optimal solution tour. \cors For instance, \corb the final step for CPP and RPP involves an algorithm to find an Eulerian cycle, see \cite{beiseltLaporte}. This property of invariance of cost with respect to the order of traversal 
might not always be applicable to all such combinatorial optimization problems; for example, problems with time-dependent travel times have different costs for different order of traversal \cors even with the same edge/arc repetitions in the solution. \corb Since RPP-TU also has time-dependent/temporal attributes, a suitable choice for decision variables is to capture the order of arc traversal. Further, the temporal attributes in RPP-TU are appropriate for modeling the unavailability of the railway tracks for pre-specified time intervals, thus capturing the effect of existing train schedules while planning the movement of an agent (inspection or maintenance vehicle). 

Regular inspection and maintenance of rail infrastructure are crucial for the reliable operation of railways. Rail track inspection along with deterioration models provides data to survey the track conditions, which, subsequently, are used for maintenance planning, as shown in \cite{xu2015}. These tasks are generally performed using inspection vehicles, see \cite{osman2018} or dispatching robots like automated visual inspection robot, presented in \cite{resendiz}. 
Most of these railway inspection and maintenance tasks are not performed during the regular operational hours of the passenger trains, thus the routing and scheduling problem under consideration doesn't have temporal properties. 
Few studies by \cite{peng}, \cite{pour}, \cite{budai} etc., model large-scale routing and scheduling problems without considering regular train schedules. However, a more versatile problem setting involves modeling the unavailability of railway tracks due to existing train schedules, for maximum utilization and reliability of the infrastructure. One such way of modeling the routing and scheduling in the presence of track unavailability, as a variant of Capacitated Arc Routing Problem (CARP), is developed in \cite{lannez}. Their approach of modeling the problem as CARP shows an increase in the number of decisions required as the temporal data (like the number of unavailability schedules) increases, and thus 
require investigation. On the other hand, the methodology employed for RPP-TU is not only useful for modeling the routing and scheduling of multiple agents without disrupting train schedules but also accounts for a fixed number of decision variables for a given problem irrespective of the temporal data. The work in \cite{b2019} describes the replicated graph and its properties given the periodically changing train schedules such that the number of decision variables is decoupled from the time limit of inspection. As a result, RPP-TU is befitting for modeling the routing and scheduling of multiple agents in a shared network, like railways. 

In order to solve time-dependent as well as time-independent combinatorial problems, they are often formulated as an Integer Programming (IP) or Mixed Integer Linear Programming (MILP) problem, which results in its solution set being represented as multiple disjoint polyhedrons. Since Linear Programming (LP) problem has been proved to be solvable in polynomial time using ellipsoidal algorithm, determining the convex hull of the feasible points to \cors formulate \corb this IP/MILP as an LP would ideally make the routing problems polynomial-time solvable. \cors For \corb example CPP is formulated as an LP using Blossom inequalities in \cite{edmonds}. However problems like TSP and RPP are known to be NP Hard. In particular, there exists no known polynomial-time solution. 
For such NP Hard problems, either deriving all the facets of the convex hull is not obvious, and/or the number of facets scales exponentially with the number of vertices; see the polyhedral study by \cite{chopra} for TSP and \cite{corberan} for RPP. Time-dependent or temporal variants of combinatorial optimization problems like Time-dependent Travelling Salesman Problem (TDTSP) by \cite{cordeau}, Time-dependent Rural Postman Problem (TDRPP) by \cite{calogiuri}, etc are NP Hard routing and scheduling problems related to roadway networks. The importance of such time-dependent traffic model for roadways is apparent in the scheduling aspect of these routing and scheduling problems, see \cite{wang2017}. Arc Path Arc Sequence (APAS) formulation by \cite{tan11a} describes a suitable methodology for representing single-agent TDRPP. Unlike these popular roadway based problem formulations, the modeling methodology for arc unavailabilities is addressed in RPP-TU which extends the APAS formulation. TDRPP is not suitable for modeling the temporal part of railway routing and scheduling problems, however, the basic structure of APAS is useful for formulating the non-temporal part of a single-agent version of RPP-TU. In addition, the polyhedral study of the feasible region of APAS formulation for TDRPP is exploited by \cite{tan11a} for improving cutting plane algorithms, like branch-and-cut. However, the extension of the polyhedral study to the multi-agent problem is non-trivial. In this work, we extend the polyhedral study by \cite{avila16} on Generalized Directed RPP to the multi-agent APAS framework. 

\noindent The main contributions of this work are enlisted as follows: 
\begin{itemize}
    \item A novel branch-and-cut algorithm is proposed that utilizes a separation algorithm to determine the cutting-planes. 
    \item Study of the polyhedral structure formed by the convex hull of the feasible solutions of the relaxation of RPP-TRU, called Cascaded Graph Formulation (CGF). The result establishes that the family of cutting-planes proposed for the branch-and-cut algorithm is facet-defining. 
    \item The criterion of the connectivity-based facet-defining inequalities is relaxed in the proposed work by allowing for the selection of disconnected sets. 
\end{itemize}

The work is structured as follows:
Section \ref{bncrpptu:sec:defs} covers a preliminary introduction to some commonly used terms and concepts. The concept of a replicated graph is also discussed in this section. Next, in Section \ref{bncrpptu:sec:cgf}, the Rural Postman Problem with Temporal Unavailability (RPP-TU) is formulated. A suitable branch-and-cut algorithm is proposed in Section \ref{bncrpptu:sec:poly_algo} to solve for an exact solution. The polyhedral study of the problem formulation is described in detail in Section \ref{bncrpptu:sec:rpptucgf} the Cascaded Graph Formulation (CGF) is also investigated in this section to prove that the proposed cutting-planes are facet-defining. 
Benchmark comparison with roadway scheduling problem and results for the simulation case study are presented in Section \ref{bncrpptu:sec:results}. Finally, the work is concluded in Section \ref{bncrpptu:sec:conc} with some propositions for future study. Detailed proofs and their illustrative examples are included in the Appendix.

%% file: sec-prelims.tex
In this section, we introduce some terminologies that are frequently referred to in our work. 

\begin{figure*}[!h]
	\centering
    \includegraphics[width=0.78\textwidth]{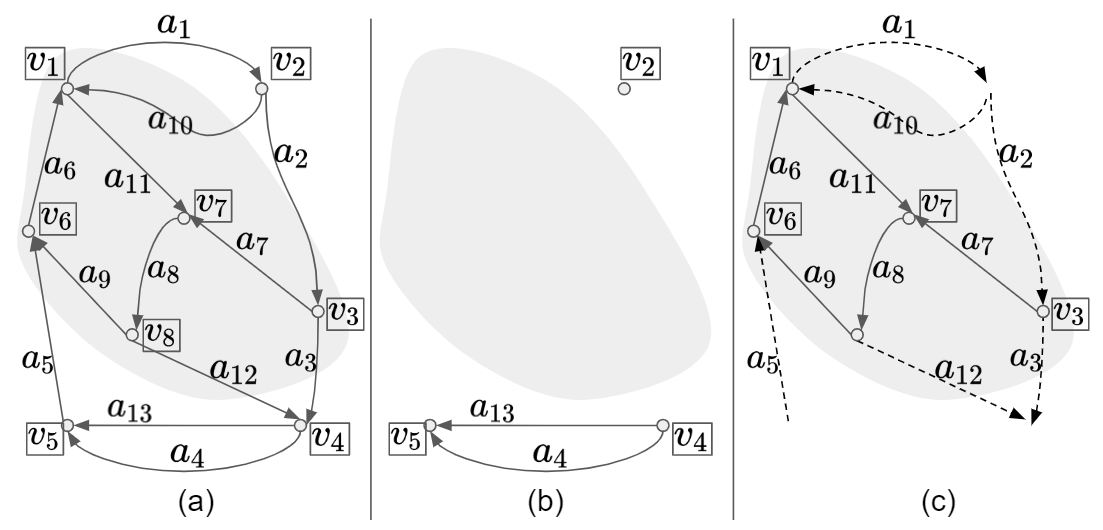}
	\captionsetup{justification=centering}
	\caption{(a) An illustrative graph $G = (V, A, F^+, F^-)$, where $V = \{v_1, \dots, v_8\}$ and $A = \{a_1, \dots, a_{13}\}$, $F^+$ and $F^-$ are suitable maps. (b) A sub-graph of $G$ with three vertices and two arcs that are strictly outside the shaded blob. (c) A sub-graph constructed by selecting all vertices and arcs strictly inside the shaded blob. The dashed arcs represent the arcs in the boundary. }
	\label{fig:prelim_fig1}
\end{figure*}

\begin{defy}{Directed multi-graph}
    A directed multi-graph is an ordered 4-tuple $G = (V, A, F^+, F^-)$; where $V$ is a vertex set, $A$ is an arc set, $F^+:A \rightarrow V$ is a map that assigns a head vertex to each arc, and $F^-:A \rightarrow V$ is a map that assigns a tail vertex to each arc. Every arc in the arc set $A$ is directed from its tail vertex to its head vertex. 
    For brevity, the terms graph and directed multi-graph are used interchangeably.  
    
    Figure \ref{fig:prelim_fig1}(a) shows a directed multi-graph with vertex set $V = \{v_1, \dots, v_8\}$ and arc set $A = \{a_1, \dots, a_{13}\}$, where $a_4$ and $a_{13}$ are parallel arcs. Both arcs $a_4$ and $a_{13}$ are  directed from their tail vertex $F^-(a_4) = F^-(a_{13}) = v_4$ to their head vertex $F^+(a_4) = F^+(a_{13}) = v_5$. 
\end{defy}

\begin{defy}{Incidence Matrix (B)}
    It is a matrix of dimension $|V|\times |A|$, that represents the relation between vertices (rows) and arcs (columns). For a directed graph, this matrix is composed of $-1$, $0$ or $1$. Every column represents an arc with one negative entry corresponding to the row at which the arc originates (-1 corresponding to the tail vertex), one positive entry corresponding to the row where the arc terminates (+1 corresponding to the head vertex), and zero for the rest. Figure \ref{fig:prelim_fig1}b shows a graph with three vertices and two arcs, and its incidence matrix is given as:
    \[
    \begin{blockarray}{ccc}
         & a_{4} & a_{13} \\
        \begin{block}{c(cc)}
            v_2 & 0 & 0 \\
            v_4 & -1 & -1 \\
            v_5 & 1 & 1 \\
        \end{block}
    \end{blockarray}
    \]
    Note that the $3 \times 2$ incidence matrix represents a multi-graph that has identical columns implying parallel arcs. 
\end{defy}

\begin{defy}{Depot}
    In a single/multi-agent problem, depot vertices indicate the physical start and end vertices in the network for the agent(s). 
    In Figure \ref{fig:prelim_fig2}, $v_1$ is considered as the depot vertex for all agents. 
    
\begin{figure}[!h]
	\centering
	\includegraphics[scale=0.35]{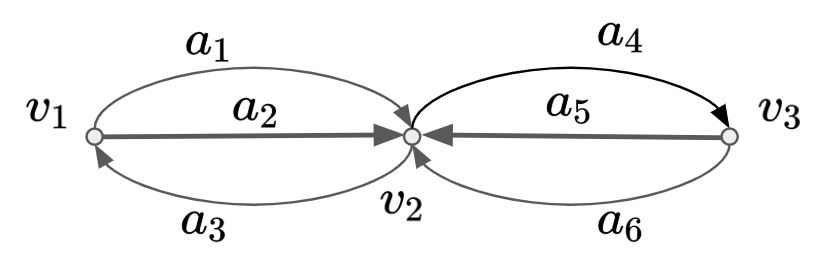}
	\captionsetup{justification=centering}
	\caption{An example graph from \cite{lannez}, where vertex $v_d = v_1$ is the depot, and $A_* := \{a_2, a_5\}$ represent the service arcs.}
	\label{fig:prelim_fig2}
\end{figure}
\end{defy}

\begin{defy}{Service arcs}
    Service arcs are a subset of the arcs set of graph $G$, denoted as $A_* \subseteq A$. A service arc represents a physical entity that requires attention, e.g. a rail-track requiring inspection. These arcs will be traversed once or multiple times, by at least one agent. 
\end{defy}

\begin{defy}{Deadhead arcs}
    Arcs that don't require servicing are called deadhead arcs, denoted as $A_D = A\backslash A_*$.
    
    The graph composed of deadhead arcs only is called \textit{deadhead sub-graph}, given as $G_D = (V, A_D, F^+, F^-)$.  
\end{defy}

\begin{defy}{Walk, Tour, Path, and Cycle (directed)}
    A (directed) \textit{walk} of length $l$ is a sequence of alternating $l+1$ vertices and $l$ arcs such that head vertex of every arc (excluding the last arc) in the sequence is a tail vertex of the next arc in the sequence. For example in Figure \ref{fig:prelim_fig1}(a), the sequence $v_1, a_1, v_2, a_{10}, v_1, a_{11}, v_7, a_8, v_8$ is a walk of length $4$. 
    
    A (directed) \textit{tour}\footnote{For brevity, a tour is allowed to have repetition of arcs, hence equivalent to a closed walk. } is a closed walk (i.e. a walk with same start and end vertex). For example in Figure \ref{fig:prelim_fig1}(a), the sequence $v_1, a_1, v_2, a_{10}, v_1, a_{11}, v_7, a_8, v_8, a_9,$ $v_6, a_6, v_1$ is a tour that starts and ends at vertex $v_1$. 
    
    A \textit{walk} without vertex repetition is called a (directed) \textit{path}. A tour without vertex repetition is called a (directed) \textit{cycle}. 
\end{defy}

\begin{defy}{Connected graph}
    A  directed graph is said to be connected if there exists a path between any two vertices. A \textit{fully connected} directed graph has a pair of opposite arcs (or single arc paths) between every pair of vertices. A \textit{strongly connected} directed graph has a pair of paths joining any two vertices in either direction.
    
    A \textit{(graph) component} is defined as a maximally connected sub-graph i.e. inclusion of any more vertex to the vertex-subset of this sub-graph will violate the connected property of this sub-graph.
    
    A graph that is not connected is called a \textit{disconnected graph}. It is a collection of all the \textit{(graph) components}. 
    Figure \ref{fig:prelim_fig1}b shows a graph with two components: one component has only one vertex $v_2$, while the other has two vertices and two arcs. 
\end{defy}


\subsubsection*{Preliminaries on replicated graph}

\cora From the literature covered in Section \ref{bncrpptu:sec:intro}, it is apparent that the order of traversal alters the cost of the solution tour in temporal/time-varying routing problems. This is resolved by representing the order of traversal, either (a) by adding another index for the decision variables along with suitable flow constraints, as observed in APAS formulation for single agent Time-Dependent Rural Postman Problem (TDRPP), introduced by \cite{tan11a}, or (b) equivalently, by introducing a graph with suitable interconnections, such that flow constraints on the decision variables defined over the arcs results in an APAS-like formulation, as observed in the case of replicated graph, first introduced by \cite{b2019}. The replicated graph is an interconnection of vertices using directed arcs, where the vertices are partitioned/categorized using layers and agent-sub-graphs. 
These layers aid in labeling repetition/copies of base graph elements, which makes the order of traversal intrinsic. Similarly, the concept of agent-sub-graph serves in separating the decisions for the agents. \corb A formal definition of replicated graph is presented below, followed by its construction from a given multi-graph, and finally categorizing all the vertices and arcs of the replicated graph based on the parent base graph.  

\begin{defy}{Replicated graph}
    A replicated graph $\mathcal{G} = (\mathcal{V}, \mathcal{A}, \mathcal{F}^+, \mathcal{F}^-)$ is a directed multi-graph, whose vertex and arc labels are carefully crafted to capture the details of a base graph $G = (V, A, F^+, F^-)$, a set of service arcs $A_* \subset A$, a depot vertex $v_d \in V$, and an agent set $\mathcal{K}$. 
\end{defy}

\vspace{1em}
\noindent\textit{Constructing a replicated graph}

    \noindent The replicated graph $\mathcal{G}$ is a processed version of base graph, that duplicates copies of the base graph elements along with introduction of few virtual vertices and arcs, based on the agent data, such that order of traversal is intrinsic for all agents.
    To construct the labels of the replicated graph, its vertices and arcs are categorized based on the agent set $\mathcal{K}$ and a layer set $\mathcal{L} := \{1,2,\dots,|A_*|+1\}$. Vertices not associated with any layer in $\mathcal{L}$ are referred to as elements of $0^{th}$ layer. 
    
\begin{figure*}[!h]
	\centering
	\includegraphics[width=0.68\textwidth]{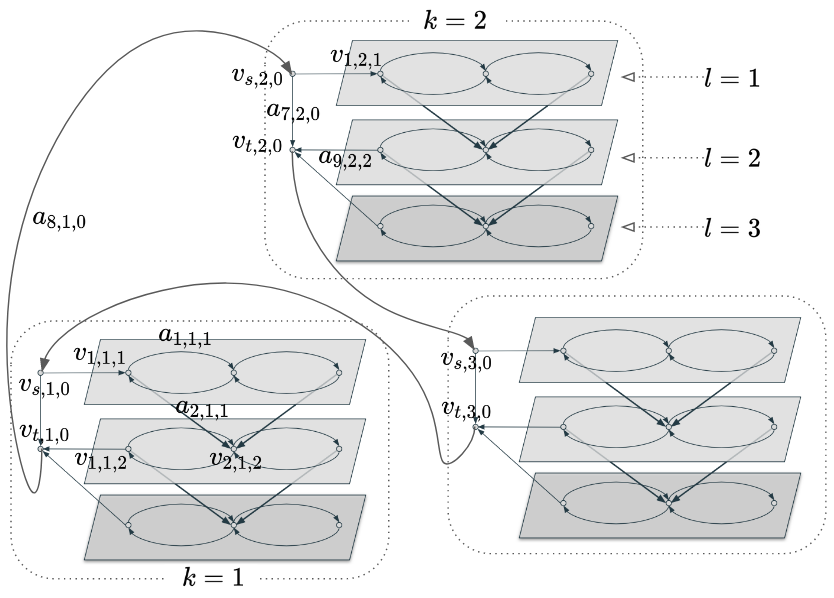}
	\captionsetup{justification=centering}
	\caption{A replicated graph for $|\mathcal{K}|$ agents, constructed from a base graph shown in Figure \ref{fig:prelim_fig2}. Here $k \in \mathcal{K} (= \{1,2,3\})$ is the agent set. The $3$ dotted boxes illustrate agent-sub-graphs having $3$ layers each \cora $l \in \mathcal{L} (= \{1,2,3\})$, \corb illustrated by planes. }
	\label{fig:prelim_fig3}
\end{figure*}

    Figure \ref{fig:prelim_fig3} shows a replicated graph with $|\mathcal{L}|$ layers (illustrated using planes) for each $|\mathcal{K}|$ \textit{agent-sub-graph}\footnote{An agent-sub-graph is a sub-graph related to an agent $k \in \mathcal{K}$.} 
    (a sub-graph, illustrated using dotted boxes). 
    Note that, in Figure \ref{fig:prelim_fig3}, the vertices are represented with three index-subscripts separated by commas, e.g. $v_{i,k,l} \in \mathcal{V}$, where $k \in \mathcal{K}$, $l \in \{0\} \cup \mathcal{L}$, and $v_i$ is either a vertex of base graph $G$ (if $l \geq 1$) or an extra vertex (if $l = 0$). For brevity, \textit{all similar notations with three index-subscripts are sometimes written without a comma separation}, e.g. $v_{ikl}$. 
    
    A pair of vertices (not belonging to the base graph vertex set $V$) called source vertex and sink vertex, are introduced for each agent-sub-graph, in the $0^{th}$ layer, resulting in $2|\mathcal{K}|$ virtual vertices, see Figure \ref{fig:prelim_fig3}. These vertices are denoted with the set $\mathcal{V}_s$ for source vertices and $\mathcal{V}_t$ for sink vertices. 
    The source and sink vertices are virtual start and end points of respective agents' route plan i.e. any agent's route starts at the source vertex $v_s \not\in V$, leads into the depot vertex $v_d \in V$ ($v_d = v_1$ in Figure \ref{fig:prelim_fig2}), performs various servicing tasks, then goes back to the depot vertex, and finally ends at the sink vertex $v_t \not\in V$. 
    
    In each layer of the replicated graph, copies of base graph vertices and deadhead arcs are embedded, resulting in $|{V}||\mathcal{K}||\mathcal{L}|$ vertices and $|{A_D}||\mathcal{K}||\mathcal{L}|$ arcs (see base graph in Figure \ref{fig:prelim_fig2} and layer illustrations in the replicated graph in Figure \ref{fig:prelim_fig3}). These vertices are called layer vertices $\mathcal{V}_l$, and the arcs are called intra-layer arcs $\mathcal{A}_l$. 
    $|\mathcal{K}|(|\mathcal{L}|-1)$ copies of the service arcs $A_*$ of the base graph $G$ are introduced between the layers such that an agent can only traverse in increasing order of the layers (top-to-bottom), see Figure \ref{fig:prelim_fig3}. These arcs are termed as inter-layer arcs $\mathcal{A}_R$, where the subscript $R$ indicates `required'. Due to this layered construction, traversing two different arcs of replicated graph might imply re-traversing the same arc of the base graph. However, in a replicated graph, the order of traversing this arc is unambiguous because one with a smaller layer index is traversed first, reflected by the top-to-bottom arc directions. 
    
    The arcs incident (incoming or outgoing) at the source and sink vertices are virtual arcs i.e. they don't represent any arc of the base graph. These virtual arcs 
    model the decision of utilizing an agents, e.g. selection of arc $a_{7,2,0}$ in Figure \ref{fig:prelim_fig3}, that connects source and sink vertices of an agent-sub-graph, implies that the $2^{nd}$ agent is at standby (no servicing jobs are assigned). 

\subsubsection*{Graph theory notations}

\begin{description}[align=right,labelwidth=0.9cm,leftmargin=1.1cm,labelsep=0.2cm,itemsep=0.0cm]
    \item [$A(S)$] indicates the set of all arcs of graph $G$ that connects the vertices in $S$. 
    Figure \ref{fig:prelim_fig1} shows an example set $S := \{v_1, v_3, v_6, v_7, v_8\}$ in the rightmost image. For this set $S$, the arcs set $A(S)$ is given as $\{a_6, a_7, a_8, a_9, a_{11}\}$. Note that $A(\mathcal{V}) = \mathcal{A}$, where $\mathcal{V}$ and $\mathcal{A}$ are vertex and arc sets of replicated graph $\mathcal{G}$.  
    \item [$G(S)$] indicates a sub-graph of $G$ with vertex set $S$ and arc set $A(S)$ i.e. $G(S) = (S, A(S), F^+, F^-)$. The sub-graph strictly inside the blob in the rightmost image of Figure \ref{fig:prelim_fig1} depicts G(S) (ignoring the dashed arcs), if $S$ is given as $\{v_1, v_3, v_6, v_7, v_8\}$. Note that $G(\mathcal{V}) = \mathcal{G}$, where $\mathcal{V}$ is a vertex set of replicated graph $\mathcal{G}$. 
    \item [$\delta(S)$] indicates the set of arcs connecting the set $S$ ($\subseteq V$) and set $V\backslash S$. In Figure \ref{fig:prelim_fig1}(a), the dashed arrows represent the boundary arcs $\delta(S) = \{a_1, a_2, a_3, a_5, a_{10}, a_{12}\}$, if $S := \{v_1, v_3, v_6, v_7, v_8\}$. 
    \item [$\delta^+(S)$] indicates the set of outgoing arcs leading from the set $S$ ($\subseteq V$) to set $V\backslash S$. In Figure \ref{fig:prelim_fig1}(a), if $S := \{v_1, v_3, v_6, v_7, v_8\}$, then $\delta^+(S) = \{a_1, a_3, a_{12}\}$. 
    \item [$\delta^-(S)$] indicates the set of incoming arcs leading from the set $V\backslash S$ to set $S$ ($\subseteq V$). In Figure \ref{fig:prelim_fig1}(a), if $S := \{v_1, v_3, v_6, v_7, v_8\}$, then $\delta^-(S) = \{a_2, a_5, a_{10}\}$. 
    \item
    \item \textit{For simplicity in notation, $\delta(v_{ikl})$ is sometimes written with a shorthand notation $\delta_{ikl}$. Similarly $\delta^+_{ikl} \implies \delta^+(v_{ikl})$ and  $\delta^-_{ikl} \implies \delta^-(v_{ikl})$.}
\end{description}

%% file: sec-rpptu.tex
The RPP-TU is a multi-agent rural postman problem with temporal restrictions due to unavailability of arcs for pre-defined periods of time. In particular, the problem description is as follows: given a graph $G=(V,A,F^+,F^-)$, running-time data for each arc $W$ \cora (arc weights), \corb a depot vertex $v_d \in V$, a set of agents $\mathcal{K}$ (such that $|\mathcal{K}| \geq 2$), an unavailability schedule $Z_q$ for each arc  $a_q \in A$ \cora (each element in the unavailability schedule set is a two tuple indicating a time interval for which the arc is inaccessible to all agents), \corb and a subset of arcs requiring service $A_* \subset A$; determine a minimum cost tour such that each service arc is traversed by at least one agent. The cost of a tour is given by the sum of fuel cost (traversed arc weights $W$) and the maximum among all finish times of the agents (see Section \ref{bncrpptu:sec:cgf}). In a railway scheduling problem, the unavailabilities listed in $Z_q$, are caused mainly due to train schedule and maintenance possession\footnote{\cors \textit{Possession} is a terminology used in railway operations and management for blocking railway track sections for some time interval, to perform inspection and/or maintenance activities. \corb }. Note that the vertices are assumed to be available, i.e. an agent may wait at any vertex for any finite period of time. A bound on the total servicing time $T_S$ is also assumed to be \cors known \corb in this work. Such a bound is usually allotted to a planner to complete all the servicing requirements for arcs in $A_*$. In practical applications, $T_S$ might be derived from the existing heuristic schedule, which is upgraded by solving RPP-TU for an exact solution. 

The RPP-TU falls in the category of time-dependent (or temporal) RPP, thus the order of traversal changes the optimal solution. 
The RPP-TU is typically formulated with two sets of variables, namely spatial and temporal variables. As the name suggests the spatial variables describe the movement (tour) of the agents in the network, while the temporal variables describe the time-stamps of this tour. 
The time-stamps are required to model the unavailability constraints at each arc, such that an agent doesn't leave a vertex while the arc ahead is unavailable for traversal. These temporal variables are also termed as departure times of the vertex; implying that the period between any two consecutive departure time of an occupied arc consists of arc traversal time followed by waiting only at the head vertex till the next arc is available. All the vertices of a feasible tour are time-stamped using temporal variables, and the tour itself is described by the spatial variables for each arc. However, a tour may repeat/re-visit many vertices, therefore a single temporal variable per vertex is not sufficient to capture the sequence/order that describes trajectory (tour + time-stamp). %

A replicated graph resolves this problem of capturing the order of traversal while re-visiting vertices, see Section \ref{bncrpptu:sec:defs}. A replicated graph $\mathcal{G} = (\mathcal{V}, \mathcal{A}, \mathcal{F}^+, \mathcal{F}^-)$ is constructed using the problem data: the base graph $G$, service arcs subset $A_*$, depot vertex $v_d$, and set of agents $\mathcal{K}$. In the replicated graph, all inter-layer arcs $\mathcal{A}_R$ are constructed using the service arcs, in a top-to-bottom fashion, to make the order of traversal intrinsic. 
The key modification in the proposed version of the replicated graph, in comparison to \cite{b2019}, is that the underlying graph is strongly connected, even for a multi-agent scenario. In particular, the graph components of replicated graph (related to agent-sub-graph of the replicated graph $\mathcal{G}$) are cascaded to form a loop i.e. sink vertex of the $k^{th}$ agent is connected to the source vertex of the $(k+1)^{th}$ agent and the sink vertex of the last agent is connected to the source vertex of the $1^{st}$ agent. 

In the following subsections, the vertex and arc subsets of replicated graph are categorized. This aids in the polyhedral study of the RPP-TRU. Next, the formulation is formally stated with brief description on the cost and constraint expressions. 

\vspace{1em}
\subsection{Categorizing vertex and arc subsets in a replicated graph}

    \noindent The replicated graph is formally categorized below using three vertex subsets $\mathcal{V} = \mathcal{V}_l \cup \mathcal{V}_s \cup \mathcal{V}_t$, and three arc subsets $\mathcal{A} = \mathcal{A}_l \cup \mathcal{A}_R \cup \mathcal{A}_F$. Mathematically, the vertex subsets are described as: 
    \begin{itemize}
        \item[1.] Layer vertices $\mathcal{V}_l := \{v_{ikl} \ | \ v_i \in V, k \in \mathcal{K}, l \in \mathcal{L}\}$ are constructed from the vertices of the base graph
        . Figure \ref{fig:prelim_fig4} shows a replicated graph with all the layer vertices placed on the planes illustrating layers. 
\begin{figure*}[!h]
	\centering
	\includegraphics[scale=0.5]{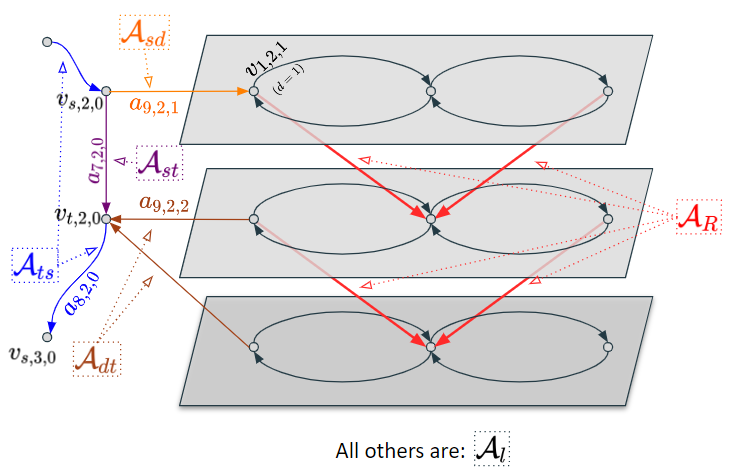}
	\captionsetup{justification=centering}
	\caption{Illustration of various categories of arcs in the sub-graph corresponding to agent $k=2$ of replicated graph $\mathcal{G}$. }
	\label{fig:prelim_fig4}
\end{figure*}
        \item[2.] Source vertices $\mathcal{V}_s := \{v_{sk0} \ | \ k \in \mathcal{K}\}$ are constructed using additional vertices $v_s \not\in V$. Figure \ref{fig:prelim_fig4} shows one source vertices for each agent-sub-graph that are contained in the $0^{th}$ layer; e.g. $v_{s,2,0} \in \mathcal{V}_s$. 
        \item[3.] Sink vertices $\mathcal{V}_t := \{v_{tk0} \ | \ k \in \mathcal{K}\}$ are constructed using additional vertices $v_t \not\in V$. Figure \ref{fig:prelim_fig4} shows one sink vertices for each agent-sub-graph that are contained in the $0^{th}$ layer; e.g. $v_{t,2,0} \in \mathcal{V}_t$. 
    \end{itemize}
    
    Similarly, the three disjoint arc subsets ($\mathcal{A}_l, \mathcal{A}_R, \mathcal{A}_F$) are mathematically described below; with further decomposition of the virtual arcs subset $\mathcal{A}_F$ into four disjoint subsets, namely source-sink arcs $\mathcal{A}_{st}$, sink-source arcs $\mathcal{A}_{ts}$, source-depodt arcs $\mathcal{A}_{sd}$, and depot-sink arcs $\mathcal{A}_{dt}$ (also shown in Figure \ref{fig:prelim_fig4}): 
    \begin{itemize}
        \item[1.] Intra-layer arcs $\mathcal{A}_l := \{a_{qkl} \in \mathcal{A} \ | \ a_q \in A_d, k \in \mathcal{K}, l \in \mathcal{L}\}$ describes all arcs that connect vertices in the same layer. It comprises all deadhead (non-service) arcs of the base graph $G$. For every $a_{qkl} \in \mathcal{A}_l$, the tail vertex is $\mathcal{F}^-(a_{qkl}) = v_{ikl}$, such that $F^-(a_q) = v_i$; and the head vertex is $\mathcal{F}^+(a_{qkl}) = v_{jkl}$, such that $F^+(a_q) = v_j$.  Figure \ref{fig:prelim_fig4} shows that these arcs are contained in their respective layers. 
        \item[2.] Inter-layer arcs $\mathcal{A}_R := \{a_{qkl} \ | \ a_q \in A_*, k \in \mathcal{K}, 1\leq l\leq |\mathcal{L}|-1\}$ comprises all arcs that connect vertices of two consecutive layers. They are constructed using service arcs $A_*$. For every $a_{qkl} \in \mathcal{A}_R$, the tail vertex is $\mathcal{F}^-(a_{qkl}) = v_{ikl}$, such that $F^-(a_q) = v_i$; and the head vertex is $\mathcal{F}^+(a_{qkl}) = v_{jkl'}$, such that $F^+(a_q) = v_j$ and $l' = l+1$. Figure \ref{fig:prelim_fig4} shows these arcs as arrows connecting one layer to the next. 
        \item[3.] Source-sink arcs $\mathcal{A}_{st}$ are given by the set $\{a_{qk0} \in \mathcal{A} \ | \ k \in \mathcal{K}\}$, for some fixed $q = \bar{q}$ such that $a_{\bar{q}} \not\in A$. These arcs connects source vertex to sink vertex, hence $\mathcal{F}^-(a_{\bar{q}k0}) = v_{sk0}$; and the head vertex is $\mathcal{F}^+(a_{\bar{q}k0}) = v_{tk0}$. Figure \ref{fig:prelim_fig4} shows that these arcs connect the source vertex subset $\mathcal{V}_s$ to the sink vertex subset $\mathcal{V}_t$; e.g. $a_{7,2,0}$ is a source-sink arc that connects $v_{s,2,0}$ to $v_{t,2,0}$ where $a_7 \not\in A$.  
        \item[4.] Sink-source arcs $\mathcal{A}_{ts}$ are given by the set $\{a_{qk0} \in \mathcal{A} \ | \ k \in \mathcal{K}\}$, for some fixed and unused $q = \tilde{q}$ ($\neq \bar{q}$) such that $a_{\tilde{q}} \not\in A$. These arcs connects sink vertex to source vertex, hence $\mathcal{F}^-(a_{\tilde{q}k0}) = v_{tk0}$; and the head vertex is $\mathcal{F}^+(a_{\tilde{q}k0}) = v_{sk'0}$, where $k' = (k+1) \ \textit{mod} \ |\mathcal{K}|$. Figure \ref{fig:prelim_fig4} shows that these arcs connect the sink vertex subset $\mathcal{V}_t$ to the source vertex subset $\mathcal{V}_s$; e.g. $a_{8,2,0}$ is a sink-source arc that connects $v_{t,2,0}$ to $v_{s,3,0}$ where $a_8 \not\in A$.
        \item[5.] Source-depot arcs $\mathcal{A}_{sd}$ are given by the set $\{a_{qk1} \in \mathcal{A} \ | \ k \in \mathcal{K}\}$, for some fixed and unused $q = \hat{q}$ ($\neq \bar{q}$, and $\neq \tilde{q}$) such that $a_{\hat{q}} \not\in A$. These arcs connect the source vertex to the depot vertex of the first layers, hence $\mathcal{F}^-(a_{\hat{q}k1}) = v_{sk0}$; and the head vertex is $\mathcal{F}^+(a_{\hat{q}k1}) = v_{dk1}$, where $v_d$ is a depot vertex. Figure \ref{fig:prelim_fig4} shows that these arcs connect the vertex subset $\mathcal{V}_s$ to the copies of depot vertex, given as $\{v_{dk1} \ | \ k \in \mathcal{K}\}$ where $d = 1$; e.g. $a_{9,2,1}$ is a source-depot arc that connects $v_{s,2,0}$ to $v_{1,2,1}$ where $a_9 \not\in A$.
        \item[6.] Depot-sink arcs $\mathcal{A}_{dt}$ are given by the set $\{a_{qkl} \in \mathcal{A} \ | \ k \in \mathcal{K}, 2\leq l\leq |\mathcal{L}|\}$, where $q = \hat{q}$ such that $a_{\hat{q}} \not\in A$. These arcs connects depot vertex to sink vertex, hence $\mathcal{F}^-(a_{\hat{q}kl}) = v_{dkl}$, where $v_d$ is a depot vertex; and the head vertex is $\mathcal{F}^+(a_{\hat{q}kl}) = v_{tk0}$. Figure \ref{fig:prelim_fig4} shows that these arcs connect the copies of depot vertex to the sink vertex subset $\mathcal{V}_t$; e.g. $a_{9,2,2}$ is a sink-depot arc that connects $v_{1,2,2}$ to $v_{t,2,0}$ where $a_9 \not\in A$. 
    \end{itemize}
    The arc set along with their respective sizes are described in Table \ref{tab:arcs_cat}. 
    \begin{table*}[h!]
        \centering
        \begin{tabular}{| p{0.75cm} | p{1.9cm} | p{1.9cm} | p{2.1cm} | p{2.1cm} | p{2.4cm} | p{2.0cm} |}
            \hline
            & Intra-layer arcs & Inter-layer arcs & Source-sink arcs & Sink-source arcs & Source-depot arcs & Depot-sink arcs\\
            \hline
            \textbf{Set} & $\mathcal{A}_l$ & \textcolor{red}{$\mathcal{A}_R$} & \textcolor{violet}{$\mathcal{A}_{st}$} & \textcolor{blue}{$\mathcal{A}_{ts}$} & \textcolor{orange}{$\mathcal{A}_{sd}$} & \textcolor{brown}{$\mathcal{A}_{dt}$} \\ 
            \hline
            \textbf{Size} & $\# 1$ & $|A_*|^2 |\mathcal{K}|$ & $|\mathcal{K}|$ & $|\mathcal{K}|$ & $|\mathcal{K}|$ & $|\mathcal{K}||A_*|$ \\ 
            \hline
        \end{tabular}
        \begin{flushleft}
            $~~~~~~~ \# 1 = |A_d| \ |\mathcal{K}| \ (|A_*|+1)$
        \end{flushleft}
        \caption{Categorizing all arcs of replicated graph}
        \label{tab:arcs_cat}
    \end{table*}
    
    Another frequently used arc subset is $\mathcal{A}_q:= \{a_{qkl} \ | \ k \in \mathcal{K}, l \in \mathcal{L}\}$ for each $a_q \in A_*$. 
    Note that, each set $\mathcal{A}_q$ contains all copies of one service arc $a_q \in A_*$, such that the union gives the required/inter-layer arc set, $\mathcal{A}_R = \cup_{a_q \in A_*} \mathcal{A}_q$. 

\subsection{\cors Three-index formulation for RPP-TU}
\label{bncrpptu:sec:cgf}

The domain set of the formulation for RPP-TU is described using the spatial variables $X \in \{0,1\}^{|\mathcal{A}|}$, and temporal variables $\Gamma \in \mathbb{R}^{|\mathcal{V}|}$; where $\mathcal{V}$ and $\mathcal{A}$ are vertex set and arc set of replicated graph $\mathcal{G}$, respectively. This results in one spatial variable for each arc of the replicated graph, and one temporal variable for each vertex of the replicated graph. 
The solution described by the spatial variables $X$ is a tour, with time stamps using temporal variables $\Gamma$ specified for all vertices included in the tour. The variables $X$ indicate whether or not an arc is occupied, while the variables $\Gamma$ indicate non-negative time values. The solution trajectory (tour + time-stamp) describes the routing plan for all agents that satisfies all the servicing requirements and network restrictions. 
In this work, we assume that the temporal attributes are absent in RPP-TU after the total servicing time $T_S$ \cora i.e. RPP-TU transforms to RPP if the start time is $T_S$. \corb Note that \cors this assumption ensures \corb that a feasible solution always exists in a time-independent APAS formulation (or replicated graph based formulation) if the underlying base graph $G$ is strongly connected.   
The existence of such a solution in a recursive unavailability scenario is discussed in the literature, see \cite{b2019}. 

The cost function captures the fuel cost\footnote{Without loss of generality, the fuel cost of any arc is assumed the same as the value of running time. } using the spatial variables, while $\gamma$ denotes the total inspection time, as described by Expression (\ref{form:costf}) and Equation (\ref{form:cost}) below. For each vertex $v_{tk0} \in \mathcal{V}_t$, the corresponding temporal variable $\Gamma_{tk0}$ is the time taken by $k^{th}$ agent ($k \in \mathcal{K}$) to complete its set of tasks and return to depot. Hence, Equation (\ref{form:cost}) ensures that $\gamma$ represents the maximum among the finishing time values of $|\mathcal{K}|$ agents. 
\begin{equation} 
\label{form:costf}
\begin{aligned}
& \underset{ }{\text{min }} & & \sum_{a_{qkl} \in \mathcal{A}_R \cup \mathcal{A}_l} W_q X_{qkl} + \gamma  
\end{aligned}
\end{equation}
\begin{equation} 
\label{form:cost}
\begin{aligned}
& \text{s.t.} & & \gamma \geq \Gamma_{tk0}, \ \ \forall v_{tk0} \in \mathcal{V}_t \\
\end{aligned}
\end{equation}
where, $W_q$ is the running-time of arc $a_q$, obtained from the vector $W$ containing all running-time data.

\noindent Now, the constraints of the three-index formulation are given as: 
\begin{itemize}
    \item[1.] Flow balance constraints: Sum of all outgoing spatial decisions (for each outgoing arc in set $\delta^+_{ikl}$) at vertex $v_{ikl}$ is same as the the sum of all incoming spatial decisions (for each incoming arc in set $\delta^-_{ikl}$) at the vertex $v_{ikl}$.
\begin{equation} 
\label{form:flow}
\begin{aligned}
 \sum_{a_{qkl} \in \delta^+_{ikl}} X_{qkl} - \sum_{a_{qkl}  \in \delta^-_{ikl}} X_{qkl} = 0,\ \ \forall v_{ikl} \in \mathcal{V} \\
\end{aligned}
\end{equation}

    \item[2.] Source constraints: Sum of the decisions of all outgoing arcs at the source vertex of $k^{th}$ agent (given by $\delta^+_{sk0}$, where $v_{sk0} \in \mathcal{V}_s$) is $1$. 
\begin{equation} 
\label{form:src}
\begin{aligned}
 \sum_{a_{qkl} \in \delta^+_{sk0}} X_{qkl} = 1,\ \ \forall v_{sk0} \in \mathcal{V}_s  \\
\end{aligned}
\end{equation}

    The source constraint ensures that there is only one integer cycle that traverses the required arcs $\mathcal{A}_R$. If there are two separate cycles in the solution, both traversing one ore more required arc, then a contradiction is observed i.e. $X_{qkl} \geq 2, ~\forall a_{qkl} \in \mathcal{A}_{ts}$ causing violation of source constraints.  

    \item[3.] Service constraints: It ensures that each of the service arc $a_q \in A_*$ is traversed only once by any one of the $|\mathcal{K}|$ agents. 
\begin{equation} 
\label{form:serv}
\begin{aligned}
 \sum_{a_{qkl} \in \mathcal{A}_q} X_{qkl} = 1,\ \ \forall a_q \in A_*  \\
\end{aligned}
\end{equation}
We denote $\mathcal{A}_q = \{a_{qkl} \in \mathcal{A}_R \ | \ k \in \mathcal{K}, l \in \mathcal{L}\}$ for each $a_q \in A_*$.

    \item[4.] Running-time constraints: It captures the temporal information from the spatial trajectories of the agents. The temporal variables $\Gamma_{jkl'}$ describe the time taken at all vertices along the path of an agent; where $l' \in \mathcal{L} \cup \{0\}$. The temporal value at any vertex can be computed by adding the running-time of an incoming arc of a solution sub-graph to the temporal value at the tail vertex of this arc. Note that, the solution sub-graph has either zero or one incoming arc, i.e. the set of occupied incoming solutions at vertex $v_{jkl'}$, denoted as $S_{jkl'} \subseteq \delta^-_{jkl'}$, is either empty or a singleton. 
\begin{equation} 
\label{form:run}
\begin{aligned}
    \Gamma_{jkl'} \geq &\tilde{W}_{q} + \Gamma_{ikl}, \\ 
    & ~~ \forall a_{qkl} \in S_{jkl'}, v_{jkl'} \not\in \mathcal{V}_s 
\end{aligned}
\end{equation}
where, 
\begin{equation} 
\nonumber
\begin{aligned}  
    S_{jkl'} &:= \{a_{qkl} \in \delta^-_{jkl'} \ | \ X_{qkl} = 1\}, \\ 
    \Gamma_{ikl} &\textit{ is the temporal value at }\\
    &\textit{  tail vertex }v_{ikl} = \mathcal{F}^-(a_{qkl}), \\
    \Gamma_{jkl'} &\textit{ is the temporal value at } \\
    &\textit{  head vertex } v_{jkl'}=\mathcal{F}^+(a_{qkl}) \\
    \tilde{W}_{q} &= \begin{cases}
    W_q &\text{$a_q \in A$}\\
    0 &\text{otherwise}
\end{cases}
\end{aligned}
\end{equation}
    Note that the source vertices, given by the set $\mathcal{V}_s$, are exempted from these set of constraints to avoid \cors conflict among the constraints. Hence, this exemption avoids infeasibility, and simultaneously \corb achieves independence in the temporal values of the agents. With this exemption, one can interpret the agent-sub-graphs to be temporally disconnected, even though the solutions are spatially connected. 



    \item[5.] Unavailability constraints: These constraints model the temporal restrictions of arcs due to existing train schedule in the network. The unavailability data is listed in $Z_q$. Any 2-tuple element $(\underline{\omega}, \overline{\omega}) \in Z_q$ describes the period in which the track section related to arc $a_q \in A$ is unavailable. For each unavailability period for arc $a_{qkl} \in \mathcal{A}_R 
    \cup \mathcal{A}_l$, constraints are modeled such that the temporal value at the tail vertex of arc $a_{qkl}$ (given by $\Gamma_{jkl}$, where $v_{jkl} = \mathcal{F}^-(a_{qkl})$) doesn't lie in the range $[\underline{\omega}-W_q, \overline{\omega}]$ if $X_{qkl} = 1$. In particular, for an occupied arc $a_{qkl}$, the departure time $\Gamma_{jkl}$ at tail vertex $v_{jkl}$ must be after the end of unavailability period ($\overline{\omega}$), or before $W_q$ units of time ahead of the arc unavailability period $\underline{\omega}$. This modification in the lower unavailability range ensures that the arc must not become unavailable within the next $W_q$ time units, while the agent is still traversing the arc.     
\begin{equation} 
\label{form:unav}
\begin{aligned}
\left.
    \begin{array}{ll}
     either,&\\
     ~~\Gamma_{jkl} \leq &\underline{\omega} - W_{q} + \tau (1 - X_{qkl})\\
     or, &\\
     ~~\Gamma_{jkl} \geq &\overline{\omega} - \tau (1 - X_{qkl}) \\
    \end{array}
\right.
\end{aligned}
\end{equation}  
\begin{equation} 
\nonumber
\begin{aligned}
\left.
    \begin{array}{l}
     \forall\; (\underline{\omega}, \overline{\omega}) \in Z_q, a_{qkl} \in \mathcal{A}_R \cup \mathcal{A}_l; \\
     \mbox{where } 0\leq l\leq |\mathcal{L}|, v_{jkl} = \mathcal{F}^-(a_{qkl})
     \end{array}
\right.
\end{aligned}
\end{equation}

\end{itemize}

Evidently, the feasible region is a polyhedron, given by the convex hull of the mixed-integer linear programming solutions of the above formulation (represented by Equations (\ref{form:cost}-\ref{form:unav})), see \cite{buriuly2022}. The \textbf{spatial sub-problem} is denoted as $F_X$, given by Equations (\ref{form:cost}-\ref{form:serv}). Later, for the polyhedral study, the spatial problem $F_X$ will be relaxed to construct CGF$_X$.

%% file: sec-algo.tex

The formulation presented in Section \ref{bncrpptu:sec:rpptucgf} models Rural Postman Problem with Temporal Unavailability (RPP-TU) with binary spatial variables and continuous (non-negative real) temporal variables. We propose a branch-and-cut algorithm that evaluates an upper and a lower bounds of an optimal solution at every node, by solving iteratively modified versions of the formulation for RPP-TU. At each node, these modifications made to the formulation are generated by reducing the feasible region using additional constraints, and then relaxing\footnote{Relaxing implies that the binary variables $X \in \{0,1\}^{|\mathcal{A}|}$ are replaced with unit boxes $0 \leq X \leq \mathbbm{1}$.} all the remaining binary variables.  

In this subsection, we first discuss the branching strategy to exhaustively obtain all possible solutions to the spatial sub-problem of the presented formulation (binary variables only). Next, we discuss the iteratively modified versions of this formulation, solved at each node of the branch-and-cut algorithm, to improve the upper and lower bound computations. A pseudo-code is provided with a step-wise explanation, along with a separation algorithm for implementing a family of valid inequalities as cutting-planes. Lastly, the polyhedral properties of the RPP-TRU are studied to establish that the proposed cutting-planes are facet-defining. 

\subsection{Branching over service arcs}
\label{bncrpptu:sec:brnch}
  
The philosophy of branching \cors is to implement \corb an exhaustive search for an optimal solution in a finite set, hence terminating at the exact solution. The spatial problem of the formulation (denoted as $F_X$) is binary programming and hence composed of a finite number of feasible solutions. 
One approach for implementing a branching strategy, in a single-agent case, is by generating one child node for each available service arc, as presented in \cite{calogiuri}. However, the branching strategy in multi-agent problems is not as direct. In this subsection, we introduce a procedure to branch over service arcs of a layered graph that addresses a multi-agent routing problem. 

\subsubsection*{Branches}

\begin{figure*}[!h]
	\centering
	\includegraphics[scale=0.55]{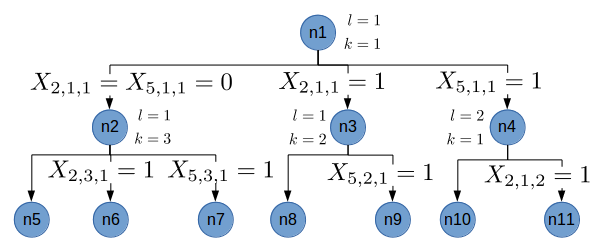}
	\captionsetup{justification=centering}
	\caption{Illustration of branching for example graph in Figure~\ref{fig:prelim_fig2}. At node $n1$, the branching is done for $l=1, k=1$; resulting in three child nodes with binary restrictions shown at the connecting branches. Assuming fractional solution is observed in $X_{2,3,1}$, the branching at node $n_2$ is performed for $l=1, k=3$. Similarly, assuming that fractional solution is observed in $X_{5,2,2}$, node $n3$ is evaluated for $l=1, k=2$ (no branching for arc $a_2$ as it has already been assigned in its parent node). If no fractional solution is achieved (say at node $n4$), then a new layer is explored for the same agent-sub-graph i.e. $l=1, k=1$.}
	\label{fig:bnc}
\end{figure*}

At each parent node, one branch is assigned for each available service arc to only one agent, and one extra branch is added to represent no assignment of service arcs to this agent (i.e. another agent will be involved next). In particular, excluding the extra branch, all other branches are assigned for each available (previously unassigned) service arc in a layer (starting from the smallest layer in $\mathcal{L}$) to only one agent (preferred based on fractional solution of $F_X$). Figure \ref{fig:bnc} shows $3$ generations of exhaustive branching for the example problem from \cite{lannez} (the example graph is shown in Figure \ref{fig:prelim_fig2}) whose replicated graph is shown in Figure \ref{fig:prelim_fig3} for $3$ agent case. Note that, in a replicated graph, if all service arc of a particular layer of agent-sub-graph are unassigned, then no succeeding layers can be assigned to this agent i.e. this agent is not available for further assignment. 

The objective of branching is to eliminate fractional solutions, hence the branching is performed on the smallest unassigned layer of any of the $|\mathcal{K}|$ agent-sub-graphs with priority for arcs with fractional spatial solution. In other words, if one of the service arcs of the smallest unassigned layers is fractional ($X_{2,3,1}$ at node $n2$ in Figure \ref{fig:bnc}) then branching is performed for this particular layer $l=1$ and agent-sub-graph $k=3$. If none of the service arcs of the smallest unassigned layer of any of the agent-sub-graphs are fractional, then either (a) the spatial solution is completely binary, and hence the branching is performed for any of the service arc of the smallest layer of an available agent-sub-graph ($n4$ in Figure \ref{fig:bnc}); or (b) there are fractional solutions in later layers ($X_{5,2,2}$ at node $n3$ in Figure \ref{fig:bnc}), hence the smallest unused layer $l=1$ of this agent-sub-graph $k=2$ is preferred for branching. 

In order to implement this branching strategy, two pieces of information are recorded at every node, (1) list of available service arcs and (2) list of available agents, for evaluating future generations of these nodes. Now, based on our branching strategy, $|A_*|+1$ child branches are generated at the root node, if the number of available agents is of size greater than $1$; otherwise there is one less child branch i.e. only $|A_*|$ branches at the root node for single-agent case. In particular, the extra branch, representing no assignment of service arcs of a layer, is not generated if the number of available agents is $1$. For each child node, if a service arc has been assigned in the previous generation of branches then the number of branches reduces by one because the list of available service arcs is shortened. If a particular child node is attached to the extra branch where no service arc of a layer was assigned, then this agent is removed from the list of available agents, and the number of child branches is same as the previous generation of branching; provided the list of available agent is of size greater than $1$. Note that, with this strategy, the number of generations in this branch-and-cut algorithm (depth of the branch-and-node search-tree) must be at most $|A_*| + |\mathcal{K}|$. 

\subsubsection*{Nodes}

The branching strategy iteratively adds nodes to a branch-and-cut search-tree. 
Some of these nodes are freshly added, and hence have no branches or child nodes in that iteration. These non-evaluated (unsolved) child nodes of the current iteration are also called \textit{leaf nodes}. 
At each iteration of the branch-and-cut algorithm, one of the leaf node with the smallest lower bound is selected and evaluated. This evaluation process involves branching to produce child nodes, and solving a modified versions of the formulation for obtaining a lower and upper bound. The modifications are due to binary relaxation, selection of service arcs in previous generations, and some improvement strategies for faster convergence to optimal solution. 


The binary relaxations convert the underlying Mixed Integer Linear Programming (MILP) problem at the node to Linear Programming (LP) problem, thus simplifying the computation. Service arc assignment to an agent is useful in reducing the number of decisions, hence further simplifying the problem to be solved at the current node. This reduction in decision variables occur because, once a service arc is selected while branching, the fastest paths to these service arcs are inherently optimal in this scenario. Hence, for each child node of the current node, a \textit{fastest path algorithm} (shortest path algorithm with running time and unavailability data, using base graph $G$) is solved, and all the arcs involved in this path are forced as $1$. All the iterative modifications are discussed in detail in Subsection \ref{bncrpptu:sec:iterF}. 
This subsection also discusses the steps for computation of an upper bound, at the current node. Note that, if the upper bound is larger than or equal to the smallest lower bound among the leaf nodes, the algorithm terminates (optimal solution is achieved). 

\subsection{Iterative modification of the formulation at the nodes}
\label{bncrpptu:sec:iterF}

Each modification of the formulation has spatial and temporal variables, which are solved separately. Observe that, given a spatial binary solution $\bar{X}$, one may substitute it into the running-time and unavailability constraints, given by Equation (\ref{form:run}) and Equation (\ref{form:unav}) respectively, to obtain constraints that depend on temporal variables $\Gamma$ only. Let $\bar{\Gamma}$ be the lowest cost temporal solution for a given binary spatial solution $\bar{X}$, then the cost of the solution $(\bar{X}, \bar{\Gamma})$ is set as the upper bound of the optimal solution. In case given $\bar{X}$ is binary but has no feasible temporal solution $\bar{\Gamma}$, or the given $\bar{X}$ is fractional, 
then the upper bound is set as the parent or the initialized upper bound. The lower bound, in all these cases, is the optimal cost of the spatial sub-problem. 
This method of solving the spatial and temporal parts separately reduces a larger problem to two smaller sub-problems. Moreover, given spatial solution $\bar{X}$, the \textit{either-or} representation of unavailability constraints, given by Equation (\ref{form:unav}) becomes non-disjoint. This simplification is observed because one of the \textit{either-or} constraints is ineffective if the spatial variable is known. 

%
The modifications in the basic formulation framework is contributed by (1) an improvement strategy that leads to better lower and upper bounds, and (2) the branching strategy described earlier that reduces the problem size by allocating decisions to some of the arcs before solving the resulting linear programming problem. 

\begin{description}[align=right,labelwidth=0.8cm,leftmargin=1.0cm, labelsep=0.2cm,itemsep=0.9cm]
    \item[\blt] The improvement strategy is incorporated into the proposed branch-and-cut algorithm to reduce the number of iterations and hence the overall computation time. 
    Cutting-planes derived from a family of proposed valid inequalities is one such strategy. The cutting-planes reduces the number of fractional solutions in the relaxation, thus reducing the number of iterations. The family of valid inequalities used as cutting-planes for $F_X$ are given in Theorem \ref{thm:conn2b_copy}. 

    
    \item[\blt] The presented branching strategy introduces another modification to our formulation, which contributes in determining a better lower bound using $|\mathcal{K}|$ additional constraints. These constraints are based on paths generated while assigning service arcs, in the branching process of the parent generations. These paths connects each agent's source vertex to various service arcs of the agents. Let these paths be denoted by $\mathcal{P}_k$, for $k \in \mathcal{K}$. Since all these path starts from the source vertex, the starting time is assumed as zero, and the total delay due to unavailabilities (difference between actual traversal time and the total running time of the paths) are computed greedily using the fastest path algorithm mentioned in the previous subsection. This time value is denoted by $d_k$, where $k \in \mathcal{K}$, and the constraints due to these paths are given by Equation (\ref{eqn:lb}). 

    \begin{equation}
        \label{eqn:lb}
        \begin{aligned}
            \hat{\gamma} \geq &\sum_{a_{qkl} \in (\mathcal{A}_R \cup \mathcal{A}_l) \backslash \mathcal{P}} W_q X_{qkl} \\
            &+ \sum_{a_{qkl} \in \mathcal{P}_k} W_q + d_k, \\
            &\textit{where, } \mathcal{P} = \cup_{k \in \mathcal{K}} \mathcal{P}_k
        \end{aligned}
    \end{equation}

    The variable $\hat{\gamma}$ acts as an estimate for the temporal cost, which produces valid lower bounds for fractional spatial solutions also. Note that, if the entire path of an agent from source to sink vertex is given, then the lower bound computed by the spatial-problem is same as the total cost (spatial cost $+$ temporal cost) of the solution.  

\end{description}

Evidently, the temporal sub-problem is expressed as:
\begin{equation}
\tag{$P_T$}
\label{form:ft}
\nonumber
\begin{aligned}
    \underset{ }{\text{min }} &\gamma \\
    \text{s.t. } &\textit{Equation (\ref{form:cost}), Equation (\ref{form:run}) and } \\ 
    &~~ \textit{Equation (\ref{form:unav})} \\ 
    & \textit{given $X$ from the spatial } \\
    &~~ \textit{ sub-problem } \\
    & \Gamma \geq 0
\end{aligned}
\end{equation}   

Similarly, the spatial-problem to be solved in each branch-and-cut iteration is of the form: 
\begin{equation} 
\tag{$P_S$}
\label{form:fx}
\nonumber
\begin{aligned}
    \underset{ }{\text{min }} &\sum_{a_{qkl} \in (\mathcal{A}_R \cup \mathcal{A}_l) \backslash \mathcal{P}} W_q X_{qkl} + \hat{\gamma} \\
    \text{s.t. } &\textit{Equations (\ref{form:flow})-(\ref{form:serv}) } \\ 
    & \textit{cutting-planes introduced iteratively } \\
    & \textit{branching constraints introduced } \\
    & ~~\textit{   iteratively } \\
    & \textit{temporal cost estimate $\hat{\gamma}$ given in } \\
    & ~~\textit{    Equation (\ref{eqn:lb}) updated iteratively } \\
    & 0 \leq X_{qkl} \leq 1, \forall a_{qkl} \in (\mathcal{A}_R \cup \mathcal{A}_l) \backslash \mathcal{P}\\
\end{aligned}
\end{equation}  

The temporal sub-problem is only evaluated if the corresponding spatial sub-problem produces binary solution and its spatial cost is smaller than the current upper bound. 

\subsection{Pseudo-code for the proposed branch-and-cut method} 
\label{bncrpptu:sec:pscode}


Algorithm \ref{algo:bnc} shows a pseudo-code to implement the proposed branch-and-cut method. Cost and constraint matrices of the presented formulation are input to the algorithm, while the output is an optimal solution. The search-tree is initialized with the input problem formulation as root node, and correspondingly, an entry of $-\infty$ is added to a queue recording lower bound of the leaf nodes, called $\textit{lbLeafNodes}$. At each node of the search-tree, data is recorded, e.g. upper bound (initialized using parent upper bound), spatial and temporal solutions (initialized as empty), forced integer variables (due to branching), parent node identity, etc. The best upper bound $ub_{best}$, and the valid inequalities are stored as common data, and updated iteratively. After adding the root node, the algorithm enters a loop which is composed of eight steps $S1-S8$ in sequence unless specified. 
\cors 
\begin{description}[align=right,labelwidth=0.7cm,leftmargin=0.9cm, labelsep=0.2cm,itemsep=0.1cm]
 \item[S1.] The first step $S1$ picks the lowest entry in $\textit{lbLeafNodes}$ which is indicative of the best lower bound so far, using function \textit{LoadNodeWithSmallestLB}. A breaking condition is also added in this step using function \textit{BreakingCondition}, which terminates the search if the duality gap between the best upper bound $ub_{best}$ and best lower bound in $\textit{lbLeafNodes}$ is zero. 
 \item[S2.] In step $S2$, the relaxation problem \ref{form:fx} is solved using function \textit{SolveLPrelaxation}, and the corresponding search-tree data is updated, including updating the corresponding entry in $\textit{lbLeafNodes}$ as $\infty$ to avoid re-selection in step $S1$. This update also ensures that the best lower bound mentioned in $S1$ is a non-decreasing function. 
 \item[S3.] Next, in step $S3$, if this lower bound is larger than the best upper bound $ub_{best}$, the algorithm returns to step $S1$ for re-selection of a better node, otherwise it jumps to next step. 
 \item[S4.] The spatial solution is checked for fractional terms in step $S4$. In case its an all integer spatial solution, the algorithm moves to step $S5$, otherwise the algorithm jumps to step $S7$ to search for valid inequalities. 
 \item[S5.] In step $S5$, a temporal solution is computed using (\ref{form:ft}), and an upper bound (the true spatio-temporal cost of the solution) is also generated, as shown by function \textit{SolveForUB}. If a temporal solution doesn't exist for this spatial solution, then its upper bound is set as the parent upper bound. 
 \item[S6.] For the feasible temporal solution scenario, if the upper bound is smaller than the best upper bound $ub_{best}$, then the $ub_{best}$ is replaced in step $S6$, otherwise it jumps to step $S8$ for branching. 
 \item[S7.] In the next step $S7$, valid inequalities are identified using function \textit{FindCuts}, also see Algorithm \ref{algo:fac}. If found, the algorithm jumps to step $S1$ to re-evaluate the relaxation problem, otherwise branch nodes are generated in step $S8$. 
 \item[S8.] This branching step $S8$ is evaluated irrespective of the spatial solution being fractional or purely integer, using function \textit{GenerateChildNodes}. Preference is given to fractional service arcs or agent-sub-graph with fractional arcs, in that order. The details of the branching strategy has been described in Subsection \ref{bncrpptu:sec:brnch}.      
\end{description} \corb 

\begin{algorithm*}
    \setstretch{1.05}
	\SetCommentSty{mycommfont}
	\SetKw{Continue}{continue}
	\SetKwFunction{solverelaxn}{SolveLPrelaxation}
	\SetKwFunction{isintgr}{IsAllInteger}
	\SetKwFunction{solvub}{SolveForUB}
	\SetKwFunction{getleafindx}{LoadNodeWithSmallestLB}
	\SetKwFunction{lookforcuts}{SeparationAlgorithm}
	\SetKwFunction{optfound}{BreakingCondition}
	\SetKwFunction{genBranches}{GenerateChildNodes}
	\SetKwFunction{validcuts}{FindCuts}
	\SetKwInOut{Input}{input}
	\SetKwInOut{InOut}{}
	\SetKwInOut{Output}{output}
	\Input{$P$ (\textit{An MILP formulation for RPP-TU}),}
	\Output{$X_{opt}, \Gamma_{opt}$}
	\tcc{Initialize}
	$T \leftarrow P$ \tcp*{add the formulation problem to search tree as root node}
	$\textit{vcuts} \leftarrow \emptyset$ \tcp*{all discovered valid cuts}
	$\textit{lbLeafNodes} \leftarrow -\infty$ \tcp*{to store lower bound of all branch-and-cut nodes}
	$ub_{best} \leftarrow \infty$ \tcp*{best upper bound}
	\tcc{Main loop}
	\While{$\textit{lbLeafNodes}$ is not entirely infinity}{
		\tcc{S1: Pick the best leaf node from $T$, and check if optimal is found}
	    $P, \textit{lbLeafNodes} \leftarrow$ \getleafindx($\textit{lbLeafNodes}$) \;
	    \optfound($lbLeafNodes, ub_{best}$) \tcp*{$ub_{best}$ is at least as large as the current best lower bound $min(lbLeafNodes)$} 
		\tcc{S2: Solve the relaxation problem and set corresponding $\textit{lbLeafNodes}$ as infinity}
	    $lb, X, T, \textit{lbLeafNodes} \leftarrow$ \solverelaxn($\textit{lbLeafNodes}$) \;
		\tcc{S3: Check if $lb$ of the current node is worth processing}
		\If{$lb \geq ub_{best}$}{
		  \Continue \tcp*{goto Step $S1$}
		}
		\tcc{S4: Check if all spatial variables in $X$ are integers}
		\eIf{\isintgr($X$)}{
		  \tcc{S5: Compute upper bound and temporal solution}
		  $ub, \Gamma \leftarrow $\solvub($P$)\;
		  \tcc{S6: Update optimal solution and $ub_{best}$}
		  \If{$ub \leq ub_{best}$}{
		    $ub_{best} \leftarrow ub$\; 
		    $X_{opt}, \Gamma_{opt} \leftarrow X, \Gamma$ 
		  }
		}{
		  \tcc{S7: Find cutting-planes and add to tree $T$}
		  $T, \textit{lbLeafNodes}, \textit{vcuts} \leftarrow$ \validcuts($T, \textit{lbLeafNodes}, \textit{vcuts}, P, X, lb$) \tcp*{see Algorithm \ref{algo:fac}}
		  \If{new valid inequalities are found}{
		    \Continue \tcp*{goto Step $S1$}
		  }
		}
		\tcc{S8: Select an agent and generate child nodes for each service arc in the smallest available layer}
	    $T, \textit{lbLeafNodes} \leftarrow$ \genBranches($P, \textit{lbLeafNodes}$) \tcp*{add child nodes in tree $T$, and corresponding leaf nodes using parent lower bound}
	}
	\caption{\textbf{crpptu}; A branch-and-cut algorithm for RPP-TU}
	\label{algo:bnc}
\end{algorithm*}

\subsubsection*{Separation algorithm} 

The number of branch exploration is reduced by adding cutting-planes to eliminate one or more fractional solutions from the LP relaxation of the current branch-and-cut node. Since this family of inequalities is {exponentially large}, a valid inequality is selected as a cutting plane only if a violation\footnote{The current fractional solution lies in the infeasible half of a cutting plane} is identified. This violation of the cutting plane is implemented using a separation algorithm. The family of valid inequalities used in this work are stated in Theorem \ref{thm:conn2b_copy}. 

\begin{thm}
    \label{thm:conn2b_copy}
    Given a service arc $a_{\breve{q}} \in A_*$ and $|\mathcal{K}| \geq 1$, the inequalities, $$\sum_{a_{qkl} \in \delta^+(S)} X_{qkl} \geq 1 - \sum_{a_{\breve{q}kl} \in H_{\breve{q}} \cap A(\mathcal{V}\backslash S)} X_{qkl}$$ are valid for the feasible spatial-polyhedron of $F_X$ that eliminate some of the fractional solutions in $F_X$, if 
    \begin{itemize}
        \item $S \subseteq \mathcal{V}\backslash\{d\}$ 
        \item $\mathcal{G} := G(\mathcal{V})$ and $G(\mathcal{V}\backslash S)$ are strongly connected
        \item $\mathcal{A}_q \cap (A(S) \cup \delta(S)) = \{\emptyset\}$ for $q \neq \breve{q}$
        \item in the sub-graph composed of arcs $A(S) \cup \delta(S)$, each component must have at least one arc from $H_{\breve{q}}$ i.e. $H_{\breve{q}} \cap (A(S_i) \cup \delta(S_i)) \neq \{\emptyset\}$.
    \end{itemize}
\end{thm}

\begin{prf}
    We show that the proposed inequalities are valid by proving that they are facet-defining in Theorem \ref{thm:conn2b}. 
    Next, the following example proves that the proposed inequalities eliminate some of the fractional solutions from RPP-TU. Given a fractional (spatial) solution for $F_X$ as $X_{\breve{q},1,1} = X_{3,1,2} = X_{\breve{q},1,2} = 0.5$ (related to path $\{a_{\breve{q},1,1}, a_{3,1,2}, a_{\breve{q},1,2}\}$ that connects two copies of service arc $a_{\breve{q}}$), then the set $S$ is constructed by selecting all the vertices composing the path. This set satisfies the properties stated in Theorem \ref{thm:conn2b_copy}, and the given fractional solution violates the inequality proposed in this Theorem. This violation is apparent from the following expressions: $\sum_{a_{qkl} \in \delta^+(S)} X_{qkl} = 0.5$ implying outgoing arc from this set has decision $0.5$, and $\sum_{a_{\breve{q}kl} \in H_{\breve{q}} \cap A(\mathcal{V}\backslash S)} X_{qkl} = 0$ (because $\sum_{a_{\breve{q}kl} \in H_{\breve{q}} \cap A(S)} X_{qkl} = 1$ implying all decisions related to arc $a_{\breve{q}kl}$ inside the set $S$ sum up to $1$). 
\end{prf}

Theorem \ref{thm:conn2b_copy} is useful for improving the computation time of our algorithm. The steps for finding this set $S$ and determining the inequality are as follows: \cors 
\begin{description}[align=right,labelwidth=1.5cm,leftmargin=1.7cm, labelsep=0.2cm,itemsep=0.1cm]
    \item[Step 1.] For every agent-sub-graph, check if more than one copy of service arc $a_{a'} \in A_*$ has a fractional solution i.e. $X_{\breve{q}kl} = \epsilon$ for more than one arc in set $H_{\breve{q}}$; for some $\epsilon \in (0, 1)$. 
    \item[Step 2.] Look for path connecting the two such arcs with fractional solution, and define its vertices as set $S$. 
    \item[Step 3.] The cutting plane is given as: sum of the arc decisions connecting the vertices of the set $S$ and the required arcs in $H_{\breve{q}}$ must be greater than or equal to $1$, i.e. $\sum_{a_{qkl} \in \delta_+(S)} X_{qkl} + \sum_{a_{qkl} \in H_{\breve{q}} \cap A(\mathcal{V}\backslash S)} X_{qkl} \geq 1$. 
\end{description} \corb 


The pseudo-code shown in Algorithm \ref{algo:fac} utilizes the \textit{SeparationAlgorithm} function to identify inequalities violated by spatial solution $X$ with respect to the problem formulation $P$. This search for a valid inequality is terminated if no such violations were found. In case a violation is found, the function \textit{UpdateTree} adds an extra node in search-tree $T$. Furthermore, the queue $\textit{lbLeafNodes}$ is updated with parent lower bound $lb$, and the valid inequalities are saved as cutting-planes in $\textit{vcuts}$. 

\begin{algorithm}
	\SetCommentSty{mycommfont}
	\SetKwFunction{updateTree}{UpdateTree}
	\SetKwFunction{lookforcuts}{SeparationAlgorithm}
	\SetKwInOut{Input}{input}
	\SetKwInOut{InOut}{}
	\SetKwInOut{Output}{output}
	\Input{$T, \textit{lbLeafNodes}, \textit{vcuts}, P, X, lb$}
	\Output{$T, \textit{lbLeafNodes}, \textit{vcuts}$}
	\BlankLine
		\While{$true$}{
		    $cuts \leftarrow$ \lookforcuts(P, X) \;
		    \eIf{$cuts$ is empty}{
		        \textit{break} \tcp*{breaks only if no more valid cuts are found}
		    }{
		        $\textit{vcuts} \leftarrow \{\textit{vcuts}, cuts\}$ \;
		        $T  \leftarrow$ \updateTree($T$)  \;
	            $\textit{lbLeafNodes} \leftarrow  \{\textit{lbLeafNodes}, lb\}$ \;
		    }
		}
	\caption{\textbf{validcuts}; Algorithm to find valid cuts}
	\label{algo:fac}
\end{algorithm}

\subsection{Polyhedral study}
\label{sec:polystudy}

In this section, the attributes of the feasible polyhedron are studied to prove that the valid inequalities in Theorem \ref{thm:conn2b_copy} are facet-defining. We first introduce a relaxation of the RPP-TRU called Cascaded Graph Formulation (CGF). Such relaxations are useful for building the theoretical base as observed in our literature survey in Section \ref{bncrpptu:sec:intro}. The relaxation only allows for multiple servicing and deadhead traversals, hence the optimal solutions of both formulations are the same. 

This study first establishes the dimension of the CGF$_X$ polyhedron in Lemma \ref{lem:dimmci}. Next, Lemma \ref{lem:conn} shows a family of facet-defining inequalities based on the dimension claim. Lastly, the Theorem \ref{thm:conn2b} shows that the family of valid inequalities proposed in Theorem \ref{thm:conn2b_copy} is facet-defining. A sketch of the proofs for the above three theoretical results is discussed in this section, while the detailed proof is included in the Appendix (see Section \ref{bncrpptu:sec:apdx}). An illustrative proof using a simple example graph is also included in the Appendix for brevity.   

\vspace{1em}
\subsubsection{Cascaded Graph Formulation (CGF)} The CGF is a relaxation of the three-index RPP-TRU formulation in Section \ref{bncrpptu:sec:cgf}. In particular, the servicing equality constraints in Equation \eqref{form:serv} are replaced with inequality constraints shown in Equation \eqref{form:serv:cgf} to allow multiple servicing, and the integer spatial constraints are relaxed to permit multiple deadhead traversals of the agents as shown in Equation \eqref{form:intlin:cgf}. 

\begin{equation}
	\label{form:serv:cgf} 
    \begin{aligned}
    	& \sum_{a_{qkl} \in \mathcal{A}_q} X_{qkl} \geq 1,\ \ \forall a_m \in A_*   \\
        &~~~~~~~\text{where, }  \mathcal{A}_q = \{a_{qkl} \in \mathcal{A}_R \ | \ k \in \mathcal{K}, l \in \mathcal{L}\}   \\
    \end{aligned}
\end{equation}

\begin{equation}
	\label{form:intlin:cgf}
    \begin{aligned}
 	& X_{qkl} \in \{0,1,2,\dots\}, \Gamma_{jkl} \geq 0 
    \end{aligned}
\end{equation}

The spatial-only formulation of CGF is expressed as:
\begin{equation}
\tag{$CGF_X$}
\label{form:cgfx}
\nonumber
\begin{aligned}
    \underset{ }{\text{min }} &\sum_{q,k,l} c_{qkl} X_{qkl} \\
    \text{s.t. } &\textit{Equations (\ref{form:flow}), (\ref{form:src}) and } \\ 
    &~~~~~~~~~~~~~~~ \textit{(\ref{form:serv:cgf}) holds} \\ 
    & X_{qkl} \in \{0,1,2,\dots\} 
\end{aligned}
\end{equation}   

The binary constraints in Expression \eqref{form:intlin:cgf} is relaxed to all non-negative integers i.e. $X_{qkl} \in \{0,1,\dots\}$. Additionally, the mathematical framework of CGF limits (upper bounds) the decision for all outgoing arcs at the source vertices by $1$ (refer to source constraints given by Equation (\ref{form:src})). 

\vspace{1em}
\subsubsection{The main theoretical results} Lemma \ref{lem:dimmci}, Lemma \ref{lem:conn} and Theorem \ref{thm:conn2b} are the three main results of the polyhedral study. 
\vspace{1em}
\begin{lem}
\label{lem:dimmci}
    For $|\mathcal{K}|\geq 2$, $dim(\text{CGF}_X) = |\mathcal{A}| - |\mathcal{V}|$. 
\end{lem}
\vspace{1em}

The dimension claim in Lemma \ref{lem:dimmci} is proved by deriving an upper and a lower bound on the dimension of the feasible polyhedron of CGF$_X$, where both bounds equal $|\mathcal{A}| - |\mathcal{V}|$. The upper bound is proved by showing $|\mathcal{V}|$ linearly independent constraints in the formulation for CGF$_X$. Hence, the upper bound on $dim(\text{CGF}_X)$ is $|\mathcal{A}| - |\mathcal{V}|$. The lower bound is proved by finding $|\mathcal{A}| - |\mathcal{V}| + 1$ affinely independent solutions, which implies there exists a subset in CGF$_X$ of dimension $|\mathcal{A}| - |\mathcal{V}|$. For counting the affinely independent solutions, the key idea adopted is to use two base solutions to ensure feasibility, and then build more solutions by modifying the base solutions. The proof is included in the Appendix, see Section \ref{bncrpptu:sec:apdx}, along with an illustrative example of the construction of a lower bound on dimension using affinely independent solutions. 
Next, a family of facet-defining inequalities is shown in Lemma \ref{lem:conn}. 

\begin{lem}
    \label{lem:conn}
    Given a service arc $a_{q'} \in A_*$ and $|\mathcal{K}| \geq 2$, the inequalities:  
    $$\sum_{a_{qkl} \in \delta^+(S)} X_{qkl} \geq 1$$ are facet inducing for CGF$_X$ polyhedron if
    \begin{itemize}
        \item $S \subseteq \mathcal{V}\backslash\{d\}$ 
        \item $\mathcal{G} := G(\mathcal{V})$ and $G(\mathcal{V}\backslash S)$ are strongly connected
        \item $\exists \mathcal{A}_{q'} \subseteq A(S) \cup \delta(S)$, and $\mathcal{A}_q \cap (A(S) \cup \delta(S)) = \{\emptyset\}$ for $m \neq m'$
        \item in the sub-graph composed of arcs $A(S) \cup \delta(S)$, each component must have at least one arc from $\mathcal{A}_{q'}$ i.e. $\mathcal{A}_{q'} \cap (A(S_i) \cup \delta(S_i)) \neq \{\emptyset\}$.
    \end{itemize}
\end{lem}

The claim in Lemma \ref{lem:conn} is proved by determining $|\mathcal{A}| - |\mathcal{V}|$ affinely independent solutions that satisfy the expression in the claim as equality. Note that set $S$ is disjoint, and the proof has to be carefully constructed considering the graphs $G(\mathcal{V}\backslash S)$, $G(S)$, and the boundary arcs $\delta(S)$. 
The detailed proof, along with an illustrative example of the construction of affinely independent solutions, is presented in the Appendix, see Section \ref{bncrpptu:sec:apdx}. 

Lastly, using the results form Lemma \ref{lem:dimmci} and Lemma \ref{lem:conn}, we show in Theorem \ref{thm:conn2b} that the vaild-inequlities proposed in Theorem \ref{thm:conn2b_copy} are facet-defining. 

\begin{thm}
    \label{thm:conn2b}
    The valid inequalities proposed in Theorem \ref{thm:conn2b_copy} are facet-defining for CGF$_X$. 
\end{thm}

\begin{prf}
    The conditions for set $S$ in this claim is same as that of Lemma \ref{lem:conn} if $\mathcal{A}_{q'} \subseteq A(S) \cup \delta(S)$ i.e. $\mathcal{A}':= \mathcal{A}_{q'} \cap A(\mathcal{V}\backslash S) = \{\emptyset\}$. For this case, $\sum_{a_{q'kl} \in \mathcal{A}'} X_{qkl} = 0$, implying that the claimed inequality $\sum_{a_{qkl} \in \delta^+(S')} X_{qkl} \geq 1 - \sum_{a_{q'kl} \in \mathcal{A}'} X_{q'kl}$ simplifies to $\sum_{a_{qkl} \in \delta^+(S)} X_{qkl} \geq 1$, which is a facet-defining inequality as per Lemma \ref{lem:conn}. 
    
    In the complementary case, when $\mathcal{A}' := \mathcal{A}_{q'} \cap A(\mathcal{V}\backslash S) \neq \{\emptyset\}$, 
    lets assume a subset $S' := \tilde{S} \cup S$ (for some arbitrary $\tilde{S}$) that satisfies the properties stated in Lemma \ref{lem:conn}. This implies $S'$ is expressed as the union of vertex subsets $S_i$ forming disconnected graphs i.e. $S' = \cup_{1 \leq i \leq r} S'_i$. The property $\mathcal{A}_{q'} \cap A(\mathcal{V}\backslash S') = \{\emptyset\}$ is also true for $S'$, thus resulting in $|\mathcal{A}| - |\mathcal{V}|$ affinely independent solutions that satisfy the equality $\sum_{a_{qkl} \in \delta^+(S')} X_{qkl} = 1 - \sum_{a_{qkl} \in \mathcal{A}_q \cap \mathcal{A}_R(\mathcal{V\backslash S'})} X_{qkl}$, where $\sum_{a_{qkl} \in \mathcal{A}_q \cap \mathcal{A}_R(\mathcal{V\backslash S'})} X_{qkl} = 0$. Observe that $\mathcal{A}'$ is contained in the set $A(\tilde{S}) \cup \delta(\tilde{S})$, because $\mathcal{A}' \subset \mathcal{A}_{q'} \subset (A(S') \cup \delta(S'))$ and $\mathcal{A}' \cap (A(S) \cup \delta(S)) = \{\emptyset\}$ (i.e. $\mathcal{A}' \subset A(\mathcal{V}\backslash S)$ as per the definition of $\mathcal{A}'$). 
    In every affinely independent solution, from Lemma \ref{lem:conn}, the following is true:  
    \begin{equation}
    \nonumber
        \begin{aligned}
            &\sum_{a_{qkl} \in \delta^+(S')} X_{qkl} = 1 \\ 
            \iff &\sum_{a_{qkl} \in \delta^+(S)} X_{qkl} + \sum_{a_{qkl} \in \delta^+(\tilde{S})} X_{qkl} = 1 \\ 
            \iff &\sum_{a_{qkl} \in \delta^+(S)} X_{qkl} = 1 - \sum_{a_{qkl} \in \delta^+(\tilde{S})} X_{qkl}
        \end{aligned}
    \end{equation}
    Recall that all the affinely independent (integer) solutions traverse only one of the required arcs in $\mathcal{A}_{q'}$ by visiting only one of the disconnected subsets $S'_i$. In particular, $\sum_{a_{qkl} \in \delta^+(\tilde{S})} X_{qkl}$ is same as $\sum_{a_{q'kl} \in \mathcal{A}'} X_{q'kl}$, which is equal to $1$, if any vertex in the set $\tilde{S}$ is visited. Hence, $\sum_{a_{qkl} \in \delta^+(S)} X_{qkl} = 1 - \sum_{a_{q'kl} \in \mathcal{A}'} X_{q'kl}$ is true for all the affinely independent solutions. The statement says: either all required arcs of the set $\mathcal{A}'$ are unoccupied resulting in $\sum_{a_{qkl} \in \delta^+(S)} X_{qkl} = 1$ and $\sum_{a_{q'kl} \in \mathcal{A}'} X_{q'kl} = 0$, or one of the required arc of the set $\mathcal{A}'$ is occupied and no arc in $\delta^+(S)$ is traversed i.e.   $\sum_{a_{qkl} \in \delta^+(S)} X_{qkl} = 0$ and $\sum_{a_{q'kl} \in \mathcal{A}'} X_{q'kl} = 1$. 
    
    
    This results in $|\mathcal{A}| - |\mathcal{V}|$ affinely independent solutions (same number as the dimension of CGF$_X$), satisfying the facet-defining expression $\sum_{a_{qkl} \in \delta^+(S)} X_{qkl} \geq 1 - \sum_{a_{q'kl} \in \mathcal{A}_{q'} \cap A(\mathcal{V\backslash S})} X_{qkl}$ as an equality.   
 \end{prf}

Lemma \ref{lem:conn} and Theorem \ref{thm:conn2b} represent the same facet of the polyhedron CGF$_X$. However, algorithmically, it is easier to select a set $S$ based on Theorem \ref{thm:conn2b}, as it is not necessary to contain any set $\mathcal{A}_q$ entirely in $A(S) \cup \delta(S)$. 

%% file: sec-res.tex
In this section, we present a comparison of our branch-and-cut algorithm for RPP-TU with branch-and-bound algorithm for TDRPP proposed by \cite{calogiuri}. This comparison involves two variants of RPP with different temporal attributes, dedicated to application in two different network types. RPP-TU is suitable for railway application with temporal attributes due to train schedules, while TDRPP is suitable for roadway application with temporal attributes due to changing traffic. Nonetheless it serves as a benchmark to assess the performance of our branch-and-cut method. We also highlight the advantage of implementing the proposed valid inequalities as cutting-planes in the branch-and-cut algorithm for RPP-TU. Once we have established a benchmark, we demonstrate the performance of our proposed branch-and-cut algorithm on a simulation case study of RPP-TU, to solve inspection routing and scheduling problem in a sub-urban railway network of Mumbai (India), without disrupting the passenger train schedules.  

\subsection{Benchmark: Comparison with TDRPP}

The proposed branch-and-cut algorithm for RPP-TU is compared with a branch-and-bound algorithm for TDRPP. The comparison justifies the performance of the proposed methodology, and hence serves as a benchmark. RPP-TU is applicable to a railway scenario while TDRPP finds application in roadways - the differences are highlighted in detail at the end of this subsection, under the heading `benchmark analysis'. 

The network parameters are adapted from the work on TDRPP by \cite{calogiuri} for fair comparison: \corh all randomly generated graphs have $20$ vertices; \corb three arc to vertex ratios, $1.2$, $1.6$, and $2$; and three service arc to deadhead arc ratio (denoted by $\beta$), $0.3$, $0.5$, and $0.7$. The number of service arcs is rounded to the nearest integer. For each of the $9$ combinations, $30$ random graphs are generated and the resulting comparison parameters are averaged. All these random temporal graph-agent data are constructed in a way that there exists a feasible solution for the problem. First, a random Hamiltonian cycle is constructed that connects all the vertices with $|V|$ arcs. Remaining arcs are connected to randomly chosen vertex pairs, and the service arcs are also randomly determined. To resemble the temporal pattern of railway network, the unavailability periods are determined only for the arcs forming Hamiltonian cycle - imitating movement of trains along this cycle. 

Table \ref{tab:compr} shows $270$ Matlab simulations of the branch-and-cut algorithm for RPP-TU (with and without the cutting-planes) and compares with that of TDRPP. The table shows data for successful instances only. \corh Simulations taking more than $1000$ seconds to converge have been marked as unsuccessful attempt, and hence skipped. \corb Benchmark data for TDRPP is directly copied from the article by \cite{calogiuri}. The comparison parameters used in the table are described as \cors follows: \corb 
\begin{description}[align=right,labelwidth=1.5cm,leftmargin=1.7cm,labelsep=0.2cm,itemsep=-0.3cm]
	\item[OPT] Number of instances (out of 30) that reached the optimal solution within a pre-specified branching limit (\textbf{$1000$ branches} - not including nodes introduced by cutting-planes), \\
	\item[NODES] Number of nodes explored in the \corh branch-and-cut \corb procedure, \\
	\item[TIME] Average computation time (in seconds) using the indicated method, \\
	\item[nCP] Average number of cutting-planes introduced.
\end{description}

{\def\arraystretch{1.13}
	\begin{table*}
		\caption{Comparison with TDRPP}
		\label{tab:compr}
		\footnotesize
		\begin{center}		
			\begin{adjustbox}{angle=0}
                \hspace{-0.8cm}
				\begin{tabular}{ccc|ccc|ccc|cccc}
\hline
& & &\multicolumn{3}{c|}{\multirow{2}{*}{TDRPP}} &\multicolumn{7}{c}{Proposed approach} \\
& & & & & &\multicolumn{3}{c}{no cutting-planes} &\multicolumn{4}{c}{with cutting-planes} \\
\hline 
$\beta$ &$|V|$ &$|A|$ &OPT &NODES &TIME (s) &OPT &NODES &TIME (s) &OPT &NODES &TIME (s) &nCP \\\hline
$0.3$ &$20$ &$24$ &$30$ &$11$ &$0.40$ &$30$ &$ 45.7$ &$ 2.26$ &$30$ &$44.90$ &$2.16$ &$0.03$ \\
& &$32$ &$30$ &$1,173$ &$30.57$ &$30$ &$ 209.8$ &$ 16.59$ &$30$ &$161.57$ &$8.959$ &$2.43$ \\
& &$40$ &$29$ &$835$ &$22.81$ &$30$ &$ 264.1$ &$ 26.11$ &$30$ &$162.40$ &$13.45$ &$3.53$ \\
$0.5$ &$20$ &$24$ &$30$ &$34$ &$0.98$ &$23$ &$174.35$ &$ 10.21$ &$30$ &$160.40$ &$12.49$ &$0.27$ \\
& &$32$ &$30$ &$2,721$ &$67.72$ &$ 3$ &$ 705$ &$ 110.25$ &$30$ &$751.27$ &$116.93$ &$8.23$ \\
& &$40$ &$22$ &$22,743$ &$614.20$ &$ 0$ &$ -$ &$ -$ &$30$ &$1462$ &$411.4$ &$28.63$ \\
$0.7$ &$20$ &$24$ &$30$ &$46$ &$1.41$ &$19$ &$552.63$ &$ 82.16$ &$30$ &$459.27$ &$60.63$ &$0.93$ \\
& &$32$ &$30$ &$11,420$ &$287.06$ &$ 0$ &$ -$ &$ -$ &$30$ &$1467.80$ &$421.78$ &$21.40$ \\
& &$40$ &$5$ &$15,285$ &$426.92$ &$ 0$ &$ -$ &$ -$ &$27$ &$3818.90$ &$2516$ &$79.48$ \\
\hline
					`$-$' & \multicolumn{11}{l}{\corh no data available for averaging \corb } \\
				\end{tabular}
			\end{adjustbox}
		\end{center}
	\end{table*}
}


\subsubsection*{Benchmark analysis} 
\label{bncrpptu:sec:disc}



RPP-TU and TDRPP are two different problems with different applications. TDRPP has properties similar to that of RPP-TU, making it an equivalently difficult problem, however there is no known way of implementing TDRPP for railway routing and scheduling problems. A close comparison of the two problems establishes a better understanding of the similarity and differences of the two problems, and consequently justify our benchmark choice. We highlight the details of this difficulty here, while discussing the assumptions made by the two problems. 

\begin{figure*}[!h]
	\centering
	\includegraphics[width=0.88\textwidth]{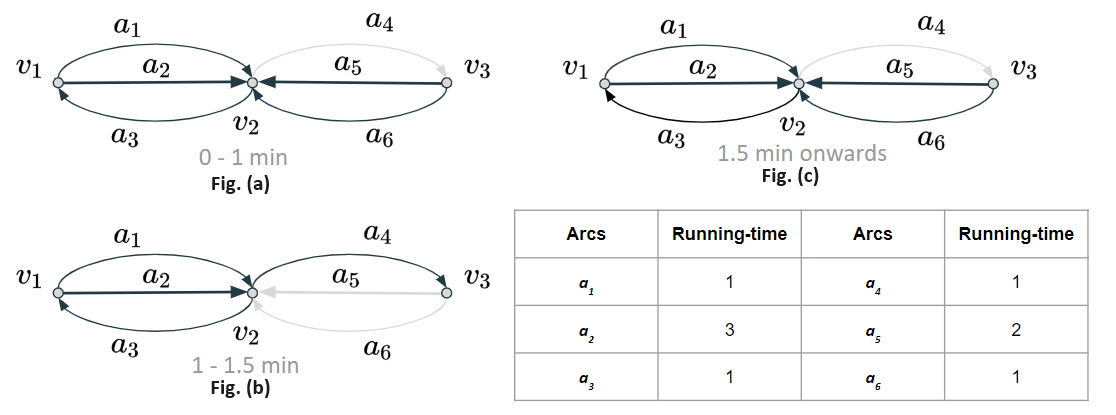}
	\captionsetup{justification=centering}
	\caption{An example showing arc $a_4$ is not available in the period $0$-$1$ minutes and also after $1.5$ minutes, while arc $a_5$ and $a_6$ are not available between $1$-$1.5$ minutes. Note that the union of the graphs is strongly connected, however arc $a_4$ is only available for $0.5$ minutes in the entire lifespan of this temporal graph. Running-time for arc $a_4$ is at least $1$ minute, hence its non-traversable. Observe that if arc $a_4$ is available in the period $0$-$1$ minutes, availability period of arc $a_4$ is larger than its running-time; however the graph is still not-traversable if $v_1$ is the depot vertex, start time is $t \geq 0$, and $a_5$ is a required arc.}
	\label{fig:rail_issues}
\end{figure*}

The roadway scheduling problems, that are modeled as TDRPP, have multiple physical factors governing its time-dependent arc-weights. Since these arc-weights represent running-time data of a road-based-agent, the weighting vector is dependent on road traffic properties like best congestion factor, degradation of this congestion factor and the running speed of the agent\footnote{The details on IGP modeling, suitable for roadway based time-dependent RPP, is proposed by \cite{ichoua}. This model assumes non-zero velocity at any point in time, which is a reasonable approximation for roadway problems.}. 
In a roadway setting, a strongly connected graph and non-zero degradation of the congestion factor is sufficient to ensure existence of a tour that satisfies all objective goals in finite time. However, in a railway setting, the underlying graph may not be connected at all time instants. Figure \ref{fig:rail_issues} shows that even though the union of the graphs at all time instants results in a strongly connected graph, availability time window of an arc might be less than running-time of the agents. Its \cors efficient \corb to attempt the modeling-optimization steps only if the problem is 
well-posed such that problems without solutions are discarded immediately. \cite{b2019} discusses some properties of railway network that serve as a base for eliminating ill-posed problems, and hence guarantee of a practically feasible solution using the RPP-TU model.    
Aside from these differences in practical setting, a first-in-first-out (FIFO) property is assumed by \cite{tan11b} and most of the succeeding work on roadway setting by \cite{tan11a}, \cite{tan12} and \cite{calogiuri}. This assumption guarantees that given two departure time from any vertex, say ${^1t_{v_1}} < {^2t_{v_1}}$ at vertex $v_1$, the time of arrival at any other vertex will be in the same order i.e. ${^1t_{v_2}} < {^2t_{v_2}}$ at vertex $v_2$. This FIFO property implies that for an optimal solution, agents need not wait at any point in the network, see \cite{tan11b}. This assumption is not reasonable in a railway setting, where trains may have different speeds, or non-uniform stop times, etc, hence resulting in optimal solution with significant amount of waiting time at various locations. 

\subsection{Simulation case study: KTVK railway network}

\begin{figure*}[!h]
	\centering
	\captionsetup{justification=centering}
	\includegraphics[width=0.88\textwidth]{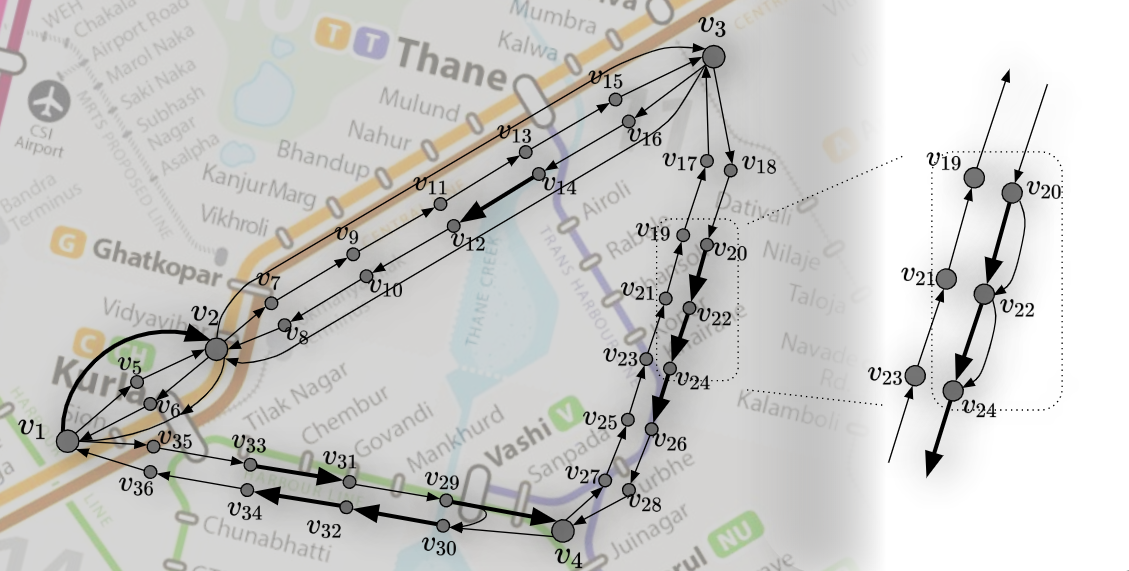}
	\caption{(This case study is borrowed from \cite{b2019}) This graph portrays 36 stations all modeled as vertices - indicated using circles of two sizes. There are 45 deadhead arcs in this model, of which 9 are to be serviced. As shown in the enlarged figure in right, the service arcs marked in bold have a duplicate deadhead arc in parallel.}
	\label{fig:ktvk}
\end{figure*}

Kurla-Thane-Vashi-Kurla (KTVK) sub-urban railway network is located in Mumbai, India. Figure \ref{fig:ktvk} shows the network map as well as a graph representation with $36$ vertices and $45$ arcs. Among them $9$ arcs that model rail-tracks require servicing (inspection). These service arcs have parallel deadhead arc representations to model motion of the agents on these tracks without performing any servicing; resulting in a total of $54$ arcs. Four types of train movements are assumed: Kurla-Thane fast and slow to-fro (between $v_1$ and $v_3$), Kurla-Vashi slow to-fro (between $v_1$ and $v_4$) and Thane-Vashi slow to-fro (between $v_1$ and $v_4$). The train schedules are assumed to be repeating every $74$ minutes. 

\noindent Table \ref{tab:table1} shows comparison results for two algorithms of RPP-TU: the proposed branch-and-cut algorithm and an earlier Benders' decomposition based algorithm (shown in \cite{b2019}). The comparison is based on three settings: 
\begin{itemize}
	\item In \textit{Ex-A}, prior knowledge of the total inspection time taken by optimal solution ($207$ minutes) is assumed as known. Hence the unavailability schedule for $222$ minutes ( = $3 \times 74$ minutes) are fed as data in the two algorithms.
	\item In \textit{Ex-B}, effect on algorithm runtime due to increase in the number of unavailabilities is observed.
	\item In \textit{Ex-C}, simulation without any unavailability scenario is presented. 
\end{itemize}

\begin{table}
	\centering
	\captionsetup{justification=centering}
	\caption{Results comparing algorithm runtime for earlier method and the proposed method. `unav' stands for unavailabilities. }
	\label{tab:table1}
    \small
	\begin{center}		
		\begin{adjustbox}{angle=0}
		\begin{threeparttable}
		\begin{tabular}{c|p{0.35in}|p{0.65in}|c|cc} 
			\hline
			& \multirow{2}{*}{No. of} & \multirow{1}{*}{Time limit } & \multirow{2}{*}{No. of} & \multicolumn{2}{c}{Algorithm runtime} \\ 
			& \multirow{2}{*}{unav} & \multirow{1}{*}{for unav} & \multirow{2}{*}{agents} & \multicolumn{2}{c}{(secs)}\\ 
			&  & (minutes) & & Earlier & Proposed \\ 
			\hline
			Ex-A & \multicolumn{1}{c|}{$719$} & 222 & 2 & 251 & 16.25 \\ 
			Ex-B & \multicolumn{1}{c|}{$4727$} & 1554 & 2 & 260 & 18.16 \\ 
			Ex-C & \multicolumn{1}{c|}{0} & - & 2 & 1.0 & 2.25 \\ 
			\hline
		\end{tabular}
		\end{threeparttable}
		\end{adjustbox}
	\end{center}
\end{table}

In both the works, the results were generated in Matlab, using an HP Probook (Core i7, 8th generation) and 8 GB of memory. Note that, in the earlier study, the temporal cost element of the optimization objective is to minimize the sum of all traversal times, however in our formulation it minimizes the max among all agents' traversal times. The key reason for the improvement observed in the computation time is because an MILP problem was solved in each iteration of the earlier method. 

The solution tour for Ex-C is same as that of a Rural Postman Problem (RPP) due to the absence of unavailability restrictions, while the solution tour for Ex-A and Ex-B are same. For Ex-A and Ex-B, the optimal tour of the first agent is: $v_1$ to $v_2$ through the service arc, $v_2$ to $v_3$ through the direct line, then to $v_4$ servicing $3$ arcs along the way, and finally to $v_1$ by servicing arc between $v_{32}$ and $v_{34}$. The optimal tour for the second agent is: $v_1$ to $v_2$ through the direct arc, $v_2$ to $v_3$ through the direct line, then back to $v_1$ (via direct line from $v_2$) and servicing an arc between $v_{12}$ and $v_{14}$ along the way, and finally servicing all the remaining arcs between $v_1$ and $v_4$. Agent-$1$ has a total runtime of $128$ minutes, where $43$ minutes is due to waiting for trains to pass by. Similarly, agent-$2$ has a total runtime of $144$ minutes, where $48$ minutes is due to waiting.

%% file: sec-apndx.tex
\section*{Proofs of Lemmas and Theorems}
\label{sec:appa}

This section is an extension of the polyhedral study for the RPP-TRU. The proof of the dimension claim in Lemma \ref{lem:dimmci}, and the facet-defining inequalities in Lemma \ref{lem:conn} is presented in complete detail. In addition, an illustrative example is included to visualize the construction of the affinely independent solutions for the proof of Lemma \ref{lem:dimmci} and Lemma \ref{lem:conn}. 

\begin{apndxprf}{Lemma \ref{lem:dimmci} (from Section \ref{sec:polystudy})}
    For $|\mathcal{K}|\geq 2$, $dim(\text{CGF}_X) = |\mathcal{A}| - |\mathcal{V}|$. 
\end{apndxprf}

\begin{prf}
    The replicated graph $\mathcal{G} := (\mathcal{V}, \mathcal{A},\mathcal{F}^+,\mathcal{F}^t)$ is a strongly connected graph. 
    This graph $\mathcal{G}$ has $|\mathcal{L}| = |A_*| + 1$ layers and $|\mathcal{K}|$ agent-sub-graphs. Note that, here $\mathcal{A}_q = \{a_{qkl} \in \mathcal{A}_R | k \in \mathcal{K}, l \in \mathcal{L}\}$ represents all copies of one particular service arc $a_q \in A_*$. As per the problem statement, every service arc must be traversed at least once, implying that one or more arcs in the $\mathcal{A}_q 
    $ must be occupied/selected for a solution to be spatially feasible. 
    
    \textbf{Proof that $|\mathcal{A}| - |\mathcal{V}|$ is an upper bound:} An upper bound on the dimension of CGF$_X$ is found by determining the number of linearly independent equations among the equality constraints. This is given by the rank of the combined coefficient matrix that models the flow constraints (i.e. the incidence matrix $\mathcal{B}$ of the replicated graph, whose rank is $|\mathcal{V}|-1$) and the source constraints ($|\mathcal{K}| \times |\mathcal{A}|$ matrix). Since the cumulative rank of the equality constraints is $|\mathcal{V}|$, an upper bound on the dimension of CGF$_X$ is $|\mathcal{A}| - |\mathcal{V}|$. %
    
    
    \textbf{Proof of $|\mathcal{A}| - |\mathcal{V}|$ is a lower bound:} A base solution tour $X_0$ is constructed for CGF$_X$ (i.e. without considering any temporal attributes) such that it occupies one required arc from each set $\mathcal{A}_q$ from the first agent-sub-graph ($k = 1$), and the source-sink arcs in all other agent-sub-graph (see $2^{nd}$ column in Table \ref{tab:lb} and Figure \ref{fig:basex0}a). In particular, $X_0$ represents a single agent solution, whose existence is guaranteed if the underlying base graph $G$ is strongly connected. For each required arc, excluding first agent-sub-graph, one affinely independent solution is constructed by connecting this arc to the tour $X_0$ without influencing the decision in the first agent-sub-graph (see $3^{rd}$ column in Table \ref{tab:lb} and compare with the $2^{nd}$ column representing base solution $X_0$; also see Figure \ref{fig:basex0}b). In order to traverse a service arc (say $a_q \in A_*$) in the $l^{th}$ layer of the second agent-sub-graph ($k = 2$), the solution $X_0$ is modified to traverse all copies of this required arcs in the preceding layers, i.e. $\{a_{q,i,2} \in \mathcal{A}_q \ | \ 1 \leq i \leq l\}$, of that agent-sub-graph (also illustrated in Figure \ref{fig:basex0}b).
    These solutions are arranged in a matrix to form an upper triangular block matrix, as shown in $3^{rd}$ column of Table \ref{tab:lb}, to give a total of $|\mathcal{A}_R| - |A_*|^2$ affinely independent solutions. 
    
    \begin{figure*}[!h]
    	\centering
    	\includegraphics[width=0.98\textwidth]{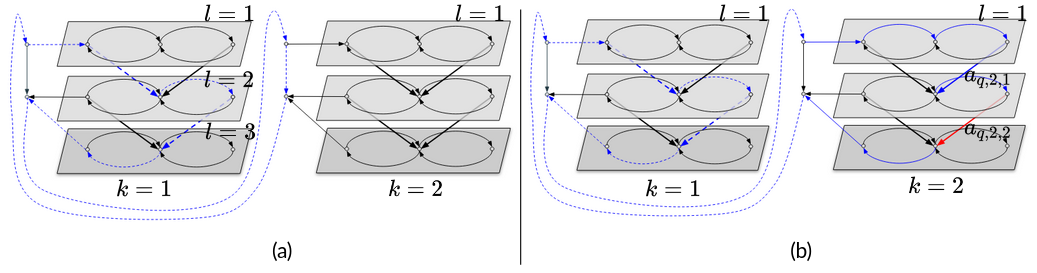}
    	\captionsetup{justification=centering}
    	\caption{(a) A base solution $X_0$, with blue-dashed arrows indicating occupied arcs. (b) Blue-dashed arrows indicate portion from the base solution $X_0$, while the blue-solid arrows indicate arcs introduced to connect the base solution to arc $a_{q,2,2} \in \mathcal{A}_q$ indicated by a red-solid arrow. Note that arc $a_{q,2,1} \in \mathcal{A}_q$, indicated by thick blue-solid arrow, is also occupied in this construction of an affinely independent solution that occupies $a_{q,2,2}$. }
    	\label{fig:basex0}
    \end{figure*} 
    
    \begin{table*}[h!]
        \centering
        \begin{tabular}{|| c || p{2.7cm} | c | c | c | c | c | c || }
            \hline
            \hline
            &  & 2 & 3 & 4 & 5 & 6 & 7s \\
            \hline
            \hline
            & & $X_0$ & using $X_0$ & $X'_{0}$ & using $X'_{0}$ & &all cycles \\ 
            \hline
            \multirow{3}{*}{2} & Required arcs of $1^{st}$ agent (occupied in $X_0$) & $\mathbbm{1}$ & $\mathbbm{1}$ & $0$ & $*$ & $*$ & $*$ \\
            \hline
            \multirow{3}{*}{3} & Depot-sink arcs of last layer for $1^{st}$ agent & \textcolor{brown}{$1$} & $1$ & $0$ & $1$ & $1$ & $1$ \\ 
            \hline
            \multirow{2}{*}{4} & Required arcs for all-but-first agent ($\mathcal{K}\backslash \{1\}$) & $0$ & \textcolor{red}{$U_{\Delta}$} & $*$ & $*$ & $*$ & $*$ \\ 
            \hline
            \multirow{2}{*}{5} & Source-sink arc for $1^{st}$ agent & $0$ & $0$ & \textcolor{violet}{$1$} & $0$ & $0$ & $*$ \\ 
            \hline
            \multirow{4}{*}{6} & Remaining required arcs for $1^{st}$ agent (unoccupied in $X_0$) & $0$ & $0$ & $0$ & \textcolor{olive}{$I_1$} & $*$ & $*$ \\
            \hline
            \multirow{3}{*}{7} & Depot-sink arcs for $1^{st}$ agent and layers $\mathcal{L}\backslash \{1,|\mathcal{L}|\}$ & $0$ & $0$ & $0$ & $0$ & \textcolor{magenta}{$I_2$} & $*$ \\ 
            \hline
            8 & $\mathcal{A}_l$ & $*$ & $*$ & $*$ & $*$ & $*$ & \textcolor{blue}{$\mathcal{C}$} \\ 
            \hline
            9 & other arcs & $*$ & $*$ & $*$ & $*$ & $*$ & $*$ \\
            \hline
            \hline
            10 & Total & $1$ & $|\mathcal{A}_R| - |A_*|^2$ & $1$ & $|A_*|^2 - |A_*|$ & $|A_*|-1$ & $\# 1$ \\
            \hline
            \hline
        \end{tabular}
        \begin{flushleft}
            \footnotesize{$\# 1$ = $(|A_D| - |V| + 1)|\mathcal{K}||\mathcal{L}|$} \\
        \end{flushleft}
        \vspace{-0.4cm}
        \captionsetup{justification=centering}
        \caption{Summarizing all affinely independent solutions in its columns. $\mathbbm{1}$ represents vector of $1$'s of size $|A_*|$, $U_{\Delta}$ represents an upper triangular matrix with all entries as $1$ or $0$; $I_1$ and $I_2$ are identity matrix of suitable size; and $\mathcal{C}$ is a matrix whose columns represent all the cycles. $*$ indicates entries not important for showing the affinely independent nature. The color code is based on the illustration in Table \ref{tab:lb_ilu}
        }
        \label{tab:lb}
    \end{table*}
    
    The next set of affinely independent solutions are constructed using the remaining required arcs that were unoccupied in the first agent-sub-graph (i.e. $k = 1$), in the solution $X_0$. To achieve this, a new base solution $X'_{0}$ is chosen that traverses all required arcs, one from each set in $\mathbb{H}$ in the second agent-sub-graph, and the source-sink arcs in all other agent-sub-graph (identical to that of $X_0$). $X'_{0}$ is affinely independent because source-sink arc of first agent is occupied for the first time (see $4^{th}$ column in the Table \ref{tab:lb}, also observe its $5^{th}$ row). The rest of the step involves constructing $|A_*|^2 - |A_*|$ affinely independent solutions, each traversing one untouched service arc from the first agent-sub-graph, but not involving the sink-source arcs $\mathcal{A}_{dt}$ of this agent-sub-graph; slightly different the earlier procedure (see $5^{th}$ column in the Table \ref{tab:lb}, also observe its $6^{th}$ and $7^{th}$ row). 

    Another $|\mathcal{L}|-2$ affinely independent solutions are constructed for each depot-sink arc $\mathcal{A}_{dt}$ in the first agent-sub-graph (except one involving last layer as shown in $7^{th}$ row of Table \ref{tab:lb}). In all the earlier solutions, these arcs were unoccupied. Hence, solutions involving these arcs adds $|\mathcal{L}| - 2 = |A_*| - 1$ affinely independent solutions to the collection, as shown in $6^{th}$ column of Table \ref{tab:lb}. 
    
    In any layer of the replicated graph, the vertices and arcs are derived from a strongly connected sub-graph $G_D ~ (\subset G)$, where $G_D := (V, A_D,F^+,F^-)$. This graph is also the base network graph with only deadhead arcs, and no service arcs. Hence, there exist $|A_D| - |V| + 1$ cycles in each layer of the replicated graph i.e. a total of $(|A_D| - |V| + 1)|\mathcal{K}||\mathcal{L}|$ cycles in the replicated graph. Since the replicated graph is strongly connected, each of these cycles can be combined with a suitable solution of CGF$_X$. Also, these cycles don't occupy any required arc, and hence cannot be created using any linear combination of the currently identified set of affinely independent solutions. 
    
    Finally we get a total of $|\mathcal{A}_R| + 1 + (|A_D| - |V| + 1)|\mathcal{K}||\mathcal{L}|$ affinely independent solutions for CGF$_X$. On simplification, the number of affinely independent solutions $n_i$ is given as:    
    \begin{equation}
        \nonumber
        \footnotesize
        \begin{aligned}
            n_i &= &&|\mathcal{A}_R|+ 1 + (|A_D| - |V| + 1)|\mathcal{K}||\mathcal{L}| \\
            &= &&(|\mathcal{A}_R| + |A_D||\mathcal{K}||\mathcal{L}| + |\mathcal{K}||\mathcal{L}|) - (|V||\mathcal{K}||\mathcal{L}|)+ 1 + 2|\mathcal{K}| - 2|\mathcal{K}| \\
            &= &&(|\mathcal{A}_R| + |A_D||\mathcal{K}||\mathcal{L}| + |\mathcal{K}||\mathcal{L}| + 2|\mathcal{K}|) - (|V||\mathcal{K}||\mathcal{L}| + 2|\mathcal{K}|) + 1 \\ 
        \end{aligned}
    \end{equation}
    Size of the arc set $|\mathcal{A}|$ is computed using Table \ref{tab:arcs_cat} in Section \ref{bncrpptu:sec:defs} ($= (|\mathcal{A}_R| + |A_D||\mathcal{K}||\mathcal{L}| + |\mathcal{K}||\mathcal{L}| + 2|\mathcal{K}|)$); and size of the vertex set $\mathcal{V}$ is the sum of the number of vertices in each layer of each agent-sub-graph ($= |V||\mathcal{K}||\mathcal{L}|$), and all the source and sink vertices ($= 2|\mathcal{K}|$). Hence, we get $n_i = |\mathcal{A}| - |\mathcal{V}|$. This gives a lower bound of $|\mathcal{A}| - |\mathcal{V}|$ on the dimension. 
\end{prf}

\begin{apndxprf}{Lemma \ref{lem:conn} (from Section \ref{sec:polystudy})} 
    Given a service arc $a_{m'} \in A_*$ and $|\mathcal{K}| \geq 2$, the inequalities:  
    $$\sum_{a_{mkl} \in \delta^+(S)} X_{mkl} \geq 1$$ are facet inducing for CGF$_X$ polyhedron if
    \begin{itemize}
        \item $S \subseteq \mathcal{V}\backslash\{d\}$ 
        \item $\mathcal{G} := G(\mathcal{V})$ and $G(\mathcal{V}\backslash S)$ are strongly connected
        \item $\exists \mathcal{A}_{m'} \subseteq A(S) \cup \delta(S)$, and $\mathcal{A}_m \cap (A(S) \cup \delta(S)) = \{\emptyset\}$ for $m \neq m'$
        \item in the sub-graph composed of arcs $A(S) \cup \delta(S)$, each component must have at least one arc from $\mathcal{A}_{m'}$ i.e. $\mathcal{A}_{m'} \cap (A(S_i) \cup \delta(S_i)) \neq \{\emptyset\}$.
    \end{itemize}
\end{apndxprf}

\begin{prf}
    Let $r$ be the number of connected components in $G(S)$ given as $\{G(S_i) | i \leq r, i \in \mathbb{N}\}$. As $\mathcal{G}$ is strongly connected, number of cut arcs $|\delta(S)|$ must be at least $2 r$. As each graph $G(S_i)$ is connected, number of arcs $|A(S_i)|$ is larger than or equal to $|S_i|-1$. 
    Since the graph $\mathcal{G}$ is strongly connected, a solution always exists that connects the set $S$ in any arbitrary way. In particular, arbitrary paths in graphs $G(S_i)$ can be joined with suitable paths in $G(\mathcal{V}\backslash S)$ to create a valid tour/solution; one set $S_i$ at a time. In each of these $r$ strongly connected graphs $G((\mathcal{V}\backslash S) \cup S_i)$ (for each $i \in \{1,\dots,r\}$), it is easy to construct a tour that satisfies all source, service and connectivity constraints, however decisions for the set $\mathcal{A}_{m'}$ must be observed carefully because it is contained in the set $A(S) \cup \delta(S)$. 
    For this proof, affinely independent solutions are constructed in two steps: first, using the sub-graph $G(\mathcal{V}\backslash S)$, with steps identical to that of Lemma \ref{lem:dimmci}; then second, using the remaining arcs $A(S) \cup \delta(S)$.
    
    \begin{table*}[!htp]
        \centering
        \begin{tabular}{|| c || p{3.7cm} | c | c | c | c || }
            \hline
            \hline
            &  & 2 & 3 & 4 & 5 \\
            \hline
            \hline
            & & $X'_{0}$ & $X_{0''}$ & using $G(\mathcal{V}\backslash S)$ & using $A(S) \cup \delta(S)$ \\ 
            \hline
            \multirow{1}{*}{2} & required arcs and cycles in $G(\mathcal{V}\backslash S)$ & $*$ & $*$ & $M_1$ & $*$ \\ 
            \hline
            \multirow{1}{*}{3} & contributing arcs in $A(S) \cup \delta(S)$ & $*$ & $*$ & $0$ & $M_2$  \\  
            \hline
            \multirow{3}{*}{4} & all other arcs of replicated graph $\mathcal{G}$ & $*$ & $*$ & $*$ & $*$ \\ 
            \hline
            \hline
            5 & Total & $-$ & $-$ & $\#1$ & $\#2$ \\ 
            \hline
            \hline
        \end{tabular}
        \begin{flushleft}
            $\# 1$ = $(|\mathcal{A}| - |A(S)| - |\delta(S)|) - (|\mathcal{V}| - |S|)$
        \end{flushleft}
        \begin{flushleft}
            $\# 2$ = $(|\delta(S)| - r) + \sum^r_{i=1} (|A(S_i)|-|S_i|+1)$
        \end{flushleft}
        \captionsetup{justification=centering}
        \caption{Counting affinely independent solutions, where $M_1$ and $M_2$ represent full rank square matrices based on the affinely independent solutions. }
        \label{tab:afsols}
    \end{table*}
    
    \textbf{Step-1}: Consider two base solutions $X_0$ and $X'_{0}$ such that $X_0$ visits the $1^{st}$ agent-sub-graph, while $X'_{0}$ visits the $2^{nd}$ agent-sub-graph, identical to that of Lemma \ref{lem:dimmci}. Without loss of generality, let $S_1$ and $S_2$ be the components visited by solutions $X_0$ and $X'_{0}$ respectively, such that the service constraints are satisfied for $\mathcal{A}_{m'}$. Now, consider the sub-graph $G(\mathcal{V} \backslash S)$ which has $|\mathcal{V} \backslash S|$ vertices and $|\mathcal{A}| - |A(S)| - |\delta(S)|$ arcs, where $|\mathcal{A}_R| - |A_*||\mathcal{K}|$ are required arcs. For each required arc in all the agent-sub-graphs except the 1st, $|\mathcal{A}_R| - |A_*||\mathcal{K}| - |A_*|(|A_*|-1)$ affinely independent solutions are constructed using base solution $X_0$ (i.e. excluding $|A_*|(|A_*|-1)$ required arcs from the $1^{st}$ agent-sub-graph); while for each required arc, that was unoccupied in $1^{st}$ agent-sub-graph based solution $X_0$, $|A_*|(|A_*|-1) - (|A_*|-1)$ affinely independent solutions are constructed using base solution $X'_{0}$ (where $|A_*|-1$ arcs were occupied in $X_0$). Furthermore, using the unoccupied depot-sink arcs of $1^{st}$ agent-sub-graph, $|A_*|-1$ affinely independent solutions are constructed. This gives a total of $|\mathcal{A}_R| - |A_*||\mathcal{K}|$ affinely independent solutions for all selection of arcs from sub-graph $G(\mathcal{V}\backslash S)$ (ignoring both $X_0$ and $X'_{0}$). 
    
    Deleting all vertices of set $S_i$ from $\mathcal{V}$ discards: (a) $(|A(S_i)|-|S_i|+1)$ cycles from $\mathcal{G}$ due to removal of arcs $A(S_i)$, and (b) $(|\delta(S_i)|-1)$ cycles more from $\mathcal{G}$ due to the removal of boundary arcs $\delta(S_i)$. This results in a total of $(|\delta(S)| - r) + \sum^r_{i=1} (|A(S_i)|-|S_i|+1)$ cycles being removed/discarded due to deleting all vertices of the set $S$. Since the required arcs are not involved in this counting of cycles, these choices must be ignored, thus resulting in removal of $(|\delta(S)| - r) + \sum^r_{i=1} (|A(S_i)|-|S_i|+1) - |A_*||\mathcal{K}|$ valid cycles from $\mathcal{G}$. Since the total number of cycles in $\mathcal{G}$ is $(|A_D| - |V| + 1)|\mathcal{K}||\mathcal{L}|$, each of the remaining $(|A_D| - |V| + 1)|\mathcal{K}||\mathcal{L}| - ((|\delta(S)| - r) + \sum^r_{i=1} (|A(S_i)|-|S_i|+1) - |A_*||\mathcal{K}|)$ cycles contribute to an affinely independent solution. Together with the contribution from required arcs, $|\mathcal{A}_R| + (|A_D| - |V| + 1)|\mathcal{K}||\mathcal{L}| - (|\delta(S)| - r) - \sum^r_{i=1} (|A(S_i)|-|S_i|+1)$ affinely independent solutions are achieved from the sub-graph $G(\mathcal{V}\backslash S)$ that satisfy the equality $\sum_{a_{mkl} \in \delta^+(S)} X_{mkl} = 1$. 
    This collection of affinely independent solutions is shown in the $4^{th}$ column of Table \ref{tab:afsols}.
    
    
    \textbf{Step-2}: On observing the cut arcs $\delta(S)$: path to each component $S_i$ and back leads to at least $r$ affinely independent solutions such that it satisfies the service constraints due to $\mathcal{A}_{m'}$ i.e. $\sum_{a_{m'kl} \in \mathcal{A}_{m'}} X_{m'kl} \geq 1$; and it visits only one of these components hence satisfying the equality $\sum_{a_{mkl} \in \delta^+(S)} X_{mkl} = 1$. Note that, in the components consisting of more than one required arc, it is possible to traverse only one of the required arcs because $G_D$ is strongly connected (implicit if the replicated graph $\mathcal{G}$ is strongly connected). The cut arcs of $\delta(S)$, that are not involved in the above step, can replace the occupied arcs of $\delta(S)$ to construct additional $|\delta(S)|-2r$ affinely independent solutions of $\mathcal{G}$ that satisfy the equality $\sum_{a_{mkl} \in \delta^+(S)} X_{mkl} = 1$. This results in a total of $|\delta(S)|-r$ affinely independent solutions that enter only one component in $G(S)$ and satisfy the given equality. 
    
    Observing components of $S$ one-by-one (say $S_i$): linearly independent directed cycles are constructed and attached to suitable paths. The number of such cycles in component $S_i$ is $|A(S_i)|-|S_i|+1$. Hence additional affinely independent solutions are constructed that satisfy the equality. 
    
    A total of $(|\delta(S)| - r) + \sum^r_{i=1} (|A(S_i)|-|S_i|+1)$ affinely independent solution are obtained using $A(S) \cup \delta(S))$, as shown in $5^{th}$ column of Table \ref{tab:afsols}. Note that, solutions $X_0$ and $X'_{0}$ are included in this collection. 
    
    On summing up the total number of affinely independent solutions $n_s$, the total results in:
    \begin{equation}
        \nonumber
        \footnotesize
        \begin{aligned}
            n_s &= &&|\mathcal{A}_R| + (|A_D| - |V| + 1)|\mathcal{K}||\mathcal{L}| \\
            &= &&(|\mathcal{A}_R| + |A_D||\mathcal{K}||\mathcal{L}| + |\mathcal{K}||\mathcal{L}|) - (|V||\mathcal{K}||\mathcal{L}|) + 2|\mathcal{K}| - 2|\mathcal{K}| \\
            &= &&(|\mathcal{A}_R| + |A_D||\mathcal{K}||\mathcal{L}| + |\mathcal{K}||\mathcal{L}| + 2|\mathcal{K}|) - (|V||\mathcal{K}||\mathcal{L}| + 2|\mathcal{K}|) 
        \end{aligned}
    \end{equation}
    Since, size of the arc set $|\mathcal{A}|$ is computed 
    while constructing the replicated graph (refer Section \ref{bncrpptu:sec:defs}) is $(|\mathcal{A}_R| + |A_D||\mathcal{K}||\mathcal{L}| + |\mathcal{K}||\mathcal{L}| + 2|\mathcal{K}|)$, and size of the vertex set $\mathcal{V}$ is the sum of the number of vertices in each layer of each agent-sub-graph ($= |V||\mathcal{K}||\mathcal{L}|$), and all the source and destination vertices ($= 2|\mathcal{K}|$), the total is $n_s = |\mathcal{A}| - |\mathcal{V}|$.  
    Hence, the total number of affinely independent solution in $\text{CGF}_X$ is $|\mathcal{A}|-|\mathcal{V}|$, therefore the inequality $\sum_{a_{mkl} \in \delta^+(S)} X_{mkl} \geq 1$ is facet defining for a suitable choice of set $S$. 
\end{prf}

Lemma \ref{lem:conn} relaxes the connectivity condition on $G(S)$, which is necessary to compensate for the fact that the connectivity of the two graphs $G(S)$ and $G(\mathcal{V}\backslash S)$, as well as the set $\mathcal{A}_{m'}$ being a subset of $A(S) \cup \delta(S)$ is not possible for multiple-agent based replicated graph. The main reason preventing such an $S$ is that the source-destination arcs can only be a part of either $G(S)$ or $G(\mathcal{V}\backslash S)$, hence only one of them is a connected graph. 





\newpage
\begin{table*}[!htp]\centering
\begin{tabular}{lp{5cm}}\toprule
\multicolumn{2}{p{6.5in}}{\textbf{Visual proof of LB for Lemma \ref{lem:dimmci}:} Through the illustrations, we show all the affinely independent solutions that contribute to the lower bound of Lemma \ref{lem:dimmci}. The illustrations also correlate the solutions with their solutions presented in Table \ref{tab:lb}. The tags in the illustrations below (e.g. R3C2) indicate the entry in the $3^\textit{rd}$ row of the $2^\textit{nd}$ column of Table \ref{tab:lb}. Observe that only the `R3C2' tag is highlighted initially, indicating that all other tagged but non-highlighted arcs have the corresponding entry as zero. These tags are get highlighted one at a time, indicating the construction of the upper triangular matrix, shown in color in Table \ref{tab:lb}. Note that, each solution presented in the illustrations below satisfies all the spatial constraints of CGF$_X$ - namely, (i) the flow is satisfied as the vertices have the same number of incoming and outgoing solid arcs, the source constraints implying there is always one solid arc leaving the source vertices, and the service constraints that ensure at least one of each category of service arcs is visited. Visually, there is one main cycle connecting the agent-sub-graphs and other small cycles are contained within a layer. } \\
\hline 
\textbf{Affinely independent solution} &\textbf{Description} \\\midrule
\multicolumn{2}{p{6.5in}}{\textbf{Row-3 Column-2 of Table \ref{tab:lb}; $\# \textit{solutions} = 1$}} \\ 
\raisebox{-\totalheight}{\includegraphics[width=0.6\textwidth, height=42mm]{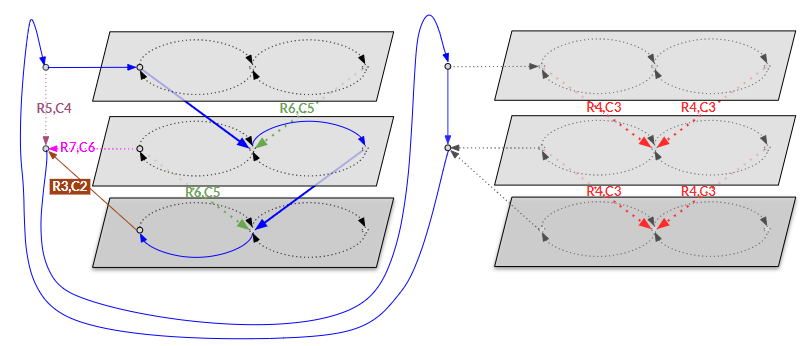}} & 
\begin{itemize}
    \setlength{\parskip}{0pt}
    \setlength{\itemsep}{0pt plus 1pt}
    \item[] Type: Base solution $X_0$ 
    \item[] Highlighted tag: R2,C2 
    \item[] Note: Observe that the solid lines that don't utilize the non-highlighted tags, hence the corresponding entries are zero. 
\end{itemize} \\
& \\ 
\multicolumn{2}{p{6.5in}}{\textbf{Row-4 Column-3 of Table \ref{tab:lb}}} \\ 
\multicolumn{2}{p{6.5in}}{$\# solutions = |\mathcal{A}_R| - |A_*||\mathcal{K}| = 4$; (\textit{i.e. Total required arcs }- \textit{Required arcs of 1st agent-sub-graph})} \\ 
\raisebox{-\totalheight}{\includegraphics[width=0.6\textwidth, height=42mm]{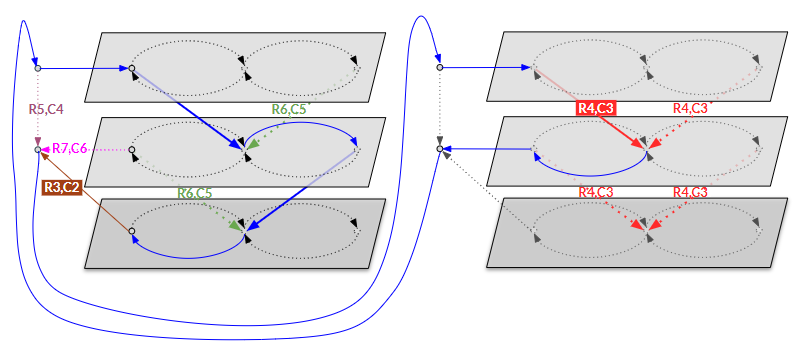}} & 
\begin{itemize}
    \setlength{\parskip}{0pt}
    \setlength{\itemsep}{0pt plus 1pt}
    \item[] Type: Using base solution $X_0$ 
    \item[] Highlighted tag: R4,C3 (1)
    \item[] Note: The modified solid blue lines now involve one additional arc (one of R4,C3). 
\end{itemize}\\ 
\raisebox{-\totalheight}{\includegraphics[width=0.6\textwidth, height=42mm]{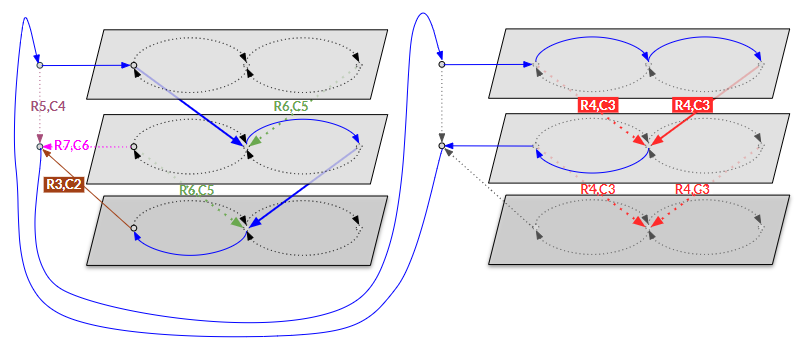}} & 
\begin{itemize}
    \setlength{\parskip}{0pt}
    \setlength{\itemsep}{0pt plus 1pt}
    \item[] Type: Using base solution $X_0$ 
    \item[] Highlighted tag: R4,C3 (2) 
    \item[] Note: The tags are checked one at a time (two of R4,C3). 
\end{itemize} \\
\raisebox{-\totalheight}{\includegraphics[width=0.6\textwidth, height=42mm]{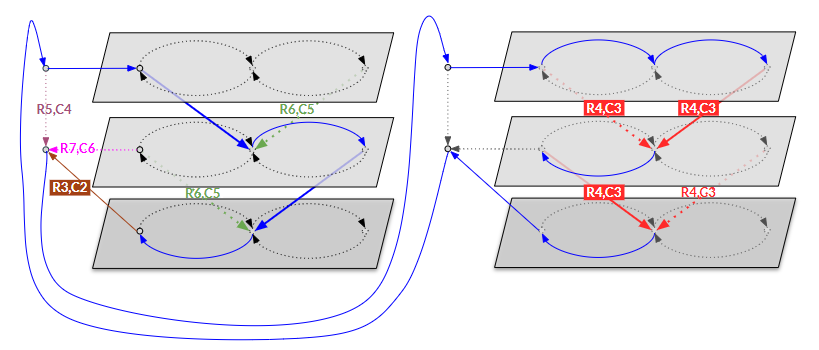}} & 
\begin{itemize}
    \setlength{\parskip}{0pt}
    \setlength{\itemsep}{0pt plus 1pt}
    \item[] Type: Using base solution $X_0$ 
    \item[] Highlighted tag: R4,C3 (3) 
    \item[] Note: We don't care whether an earlier tagged arc is revisited or not.
\end{itemize} \\
\bottomrule
\end{tabular}
\end{table*}

\newpage
\begin{table*}[!htp]\centering
\begin{tabular}{lp{5cm}}\toprule
\textbf{Affinely independent solution} &\textbf{Description} \\\midrule
\raisebox{-\totalheight}{\includegraphics[width=0.6\textwidth, height=42mm]{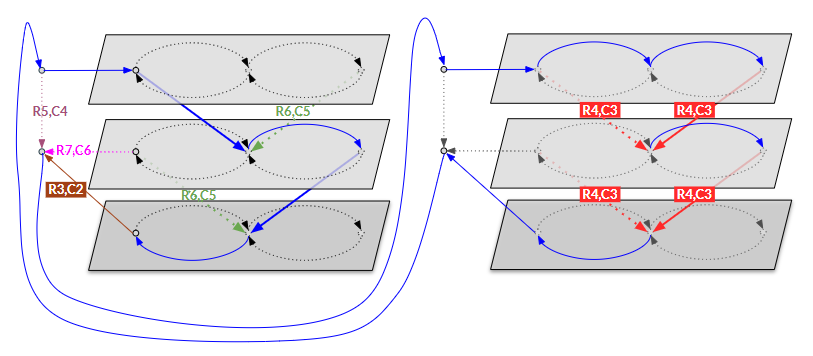}} & 
\begin{itemize}
    \setlength{\parskip}{0pt}
    \setlength{\itemsep}{0pt plus 1pt}
    \item[] Type: Using base solution $X_0$ 
    \item[] Highlighted tag: R4,C3 (4)
\end{itemize} \\
\multicolumn{2}{p{6.5in}}{\textbf{Row-5 Column-4 of Table \ref{tab:lb}; $\#\textit{solutions} = 1$}} \\ 
\raisebox{-\totalheight}{\includegraphics[width=0.6\textwidth, height=42mm]{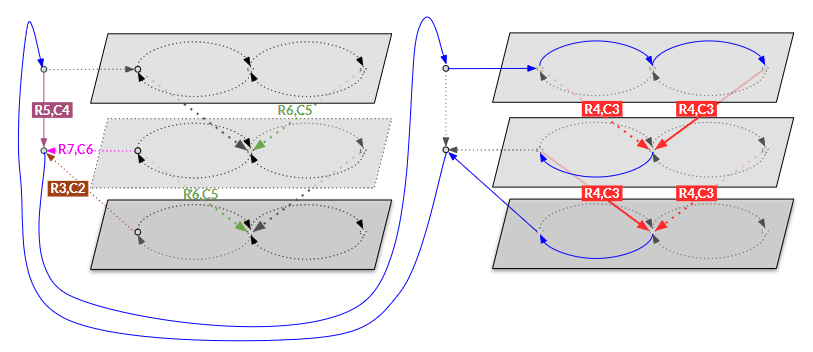}} & 
\begin{itemize}
    \setlength{\parskip}{0pt}
    \setlength{\itemsep}{0pt plus 1pt}
    \item[] Type: Base solution $X_0'$
    \item[] Highlighted tag: R5,C4 
    \item[] Note: None of the arcs in the $1^\textit{st}$ agent-sub-graph contribute to the solution. 
\end{itemize} \\
& \\ 
\multicolumn{2}{p{6.5in}}{\textbf{Row-6 Column-5 of Table \ref{tab:lb}}} \\ 
\multicolumn{2}{p{6.5in}}{$\# solutions = |A_*|(|A_*| - 1) = 2$; (\textit{i.e. Required arcs of 1st agent-sub-graph }- \textit{Required arcs already visited earlier by $X_0$})} \\ 
\raisebox{-\totalheight}{\includegraphics[width=0.6\textwidth, height=42mm]{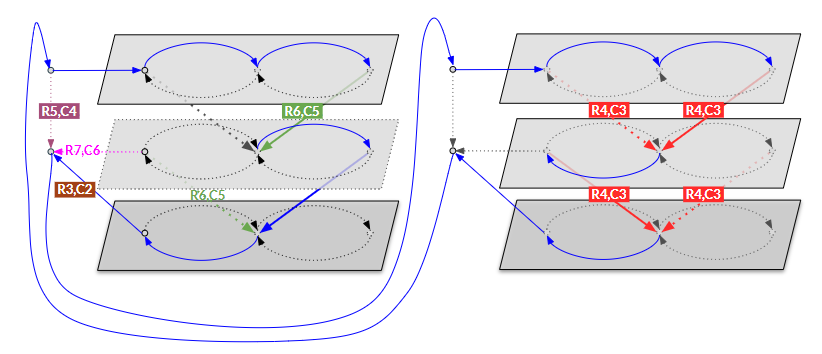}} & 
\begin{itemize}
    \setlength{\parskip}{0pt}
    \setlength{\itemsep}{0pt plus 1pt}
    \item[] Type: Using base solution $X_0'$ 
    \item[] Highlighted tag: R6,C3 (1) 
    \item[] Note: The remaining required arcs of the $1^\textit{st}$ agent-sub-graph are tagged. 
\end{itemize} \\
\raisebox{-\totalheight}{\includegraphics[width=0.6\textwidth, height=42mm]{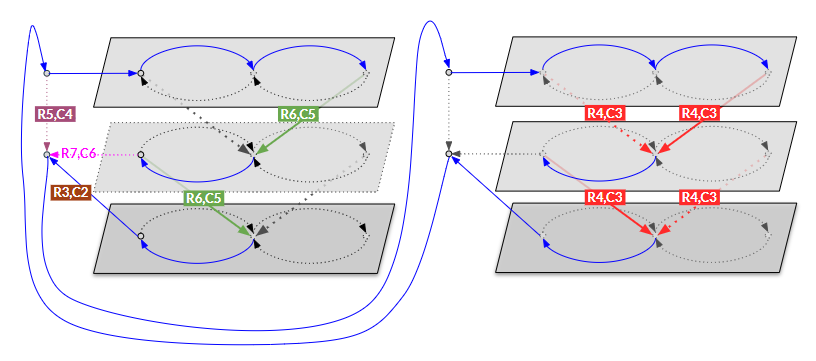}} & 
\begin{itemize}
    \setlength{\parskip}{0pt}
    \setlength{\itemsep}{0pt plus 1pt}
    \item[] Type: Using base solution $X_0'$ 
    \item[] Highlighted tag: R6,C3 (2) 
\end{itemize} \\
\multicolumn{2}{p{6.5in}}{\textbf{Row-7 Column-6 of Table \ref{tab:lb}; $\# solutions = |A_*|-1 = 1$, (for each depot-sink arc)}} \\ 
\raisebox{-\totalheight}{\includegraphics[width=0.6\textwidth, height=42mm]{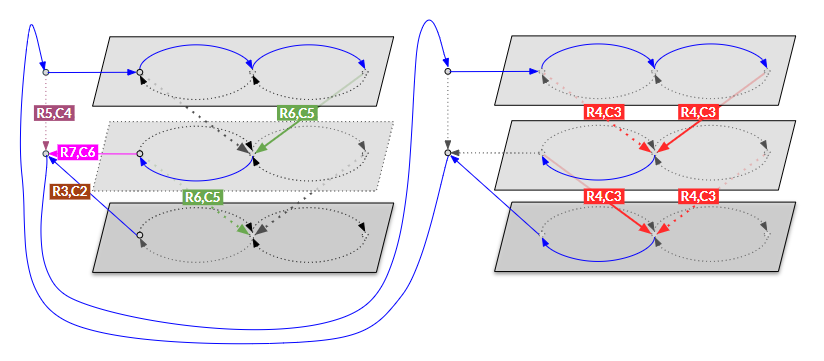}} & 
\begin{itemize}
    \setlength{\parskip}{0pt}
    \setlength{\itemsep}{0pt plus 1pt}
    \item[] Type: Using base solution $X_0'$ 
    \item[] Highlighted tag: R6,C3 (2) 
\end{itemize} \\
\bottomrule
\end{tabular}
\end{table*}

\newpage
\begin{table*}[!htp]\centering
\begin{tabular}{lp{5cm}}\toprule
\textbf{Affinely independent solution} &\textbf{Description} \\\midrule
\multicolumn{2}{p{6.5in}}{\textbf{Row-8 Column-7 of Table \ref{tab:lb} = $\# solutions = (|A_D| - |V| + 1)|A_*||\mathcal{K}| = 12$ (all intra-layer cycles)}} \\ 
\raisebox{-\totalheight}{\includegraphics[width=0.6\textwidth, height=42mm]{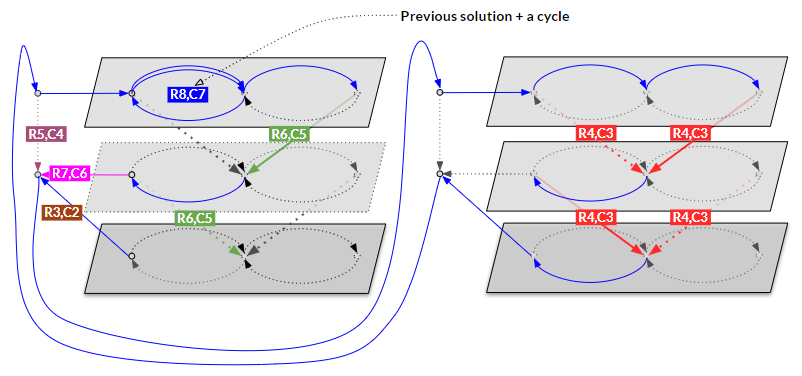}} & 
\begin{itemize}
    \setlength{\parskip}{0pt}
    \setlength{\itemsep}{0pt plus 1pt}
    \item[] Type: Using base solution $X_0$ and $X_0'$ 
    \item[] Highlighted tag: R8,C7 (1) 
\end{itemize} \\
\raisebox{-\totalheight}{\includegraphics[width=0.6\textwidth, height=42mm]{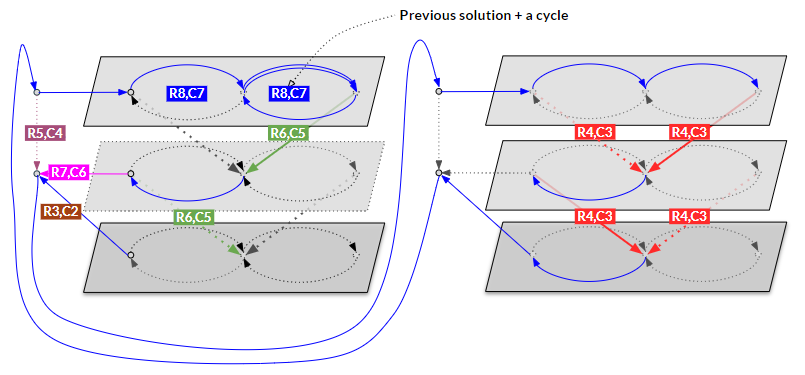}} & 
\begin{itemize}
    \setlength{\parskip}{0pt}
    \setlength{\itemsep}{0pt plus 1pt}
    \item[] Type: Using base solution $X_0$ and $X_0'$ 
    \item[] Highlighted tag: R8,C7 (2)  
\end{itemize} \\
& \\ 
\multicolumn{2}{p{6.5in}}{\textbf{$\vdots$}} \\ 
\raisebox{-\totalheight}{\includegraphics[width=0.6\textwidth, height=42mm]{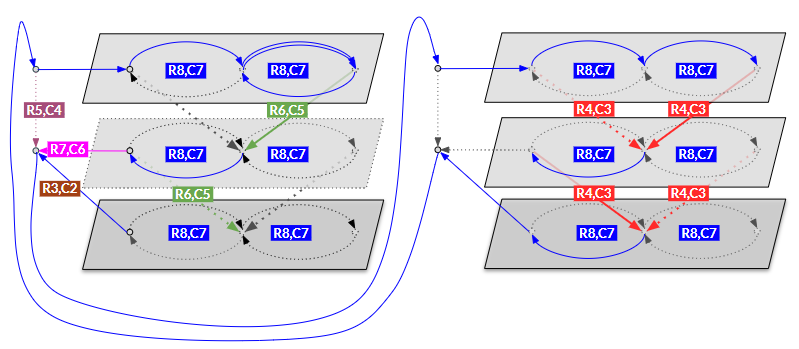}} & 
\begin{itemize}
    \setlength{\parskip}{0pt}
    \setlength{\itemsep}{0pt plus 1pt}
    \item[] Type: Using base solution $X_0$ and $X_0'$ 
    \item[] Highlighted tag: R8,C7 (12)  
\end{itemize} \\
\hline
\multicolumn{2}{p{6.5in}}{\textbf{Total number of affinely independent solution = $21$, therefore $\textit{dim$(\text{CGF}_X)$} \geq 20$}} \\ 
\bottomrule
\end{tabular}
\end{table*}


\newpage
\begin{table*}[h!]
    \centering
    \captionsetup{justification=centering}
    \caption{\textbf{For proving facet-defining cuts from Lemma \ref{lem:conn}:} Summarizing all affinely independent solutions in its columns. $U_{\Delta}^j$ represents an upper triangular matrix with all entries as $1$ or $0$, and $\mathcal{C}_1 \& \mathcal{C}_2$ are matrices whose columns represent all the cycles. $*$ indicates entries not important for showing the affinely independent nature. The color code is based on the illustrations below. In Table \ref{tab:afsols}, the matrix M1 is constructed by columns $3,5,7, \& 9$ from below, while the matrix M2 is constructed by columns $2,4,6,8, \& 10$ from below.}
    \label{tab:lbafsol}
    \begin{tabular}{|| c || p{1.0cm} | p{2.1cm} | c | c | c | c | c | c | c | c | c || }
        \hline
        \hline
        & & &2 &3 &4 &5 &6 &7 &8 &9 &10 \\
        \hline
        \hline
        & & &\multirow{2}{*}{X0} &using &X0' and &using &using &using &using &\multirow{2}{*}{cycles} &\multirow{2}{*}{cycles} \\
        & & & &X0 &using X0' &X0' &X0' &X0' &X0 or X0' & & \\
        \hline
        2 &Step-2 &Depot-sink arc of last layer of 1st agent-sub-graph & \textcolor{brown}{1} &0 &* &* &* &* &* &* &* \\
        \hline
        3 &Step-1 &Required arcs in k-th agent-sub-graph contained in G(V\S); (k >= 2) &0 &\textcolor{red}{$U_\Delta^1$} &* &* &* &* &* &* &* \\
        \hline
        4 &Step-2 &Required arcs in k-th agent-sub-graph utilizing Si or boundary; (k>=2) &0 &0 &\textcolor{red}{$U_\Delta^2$} &* &* &* &* &* &* \\
        \hline
        5 &Step-1 &Required arcs in 1st agent-sub-graph contained in G(V\S) &0 &0 &0 &\textcolor{black!30!green}{$U_\Delta^3$} &* &* &* &* &* \\
        \hline
        6 &Step-2 &Required arcs in 1st agent-sub-graph utilizing Si or boundary &0 &0 &0 &0 &\textcolor{black!30!green}{$U_\Delta^4$} &* &* &* &* \\
        \hline
        7 &Step-1 &Depot-sink arcs in 1st agent subgraph &0 &0 &0 &0 &0 &\textcolor{red!30!pink}{$U_\Delta^5$} &* &* &* \\
        8 &Step-2 &Other boundary arcs &0 &0 &0 &0 &0 &0 &$U_\Delta^6$ &* &* \\
        \hline
        9 &Step-1 (cyc) &Cycles in G(V\S) &* &* &* &* &* &* &* &\textcolor{gray}{C1} &* \\
        \hline
        10 &Step-2 (cyc) &Cycles in G(S) &* &* &* &* &* &* &* &* &\textcolor{gray}{C2} \\
        \hline
        11 & &Rest of the arcs &* &* &* &* &* &* &* &* &* \\
        \hline
        \hline
        12 &\multicolumn{2}{|c|}{\textbf{Total}} &1 &$\# 1$ &$r-|A_*|$ & $\# 2$ &$|A_*|-1$ &$|A_*|-1$ &$|\delta(S)|-2r$ & $\# 3$ & $\# 4$ \\
        \hline
        \hline
        13 &\multicolumn{2}{|c|}{Total (in example) = 20} &1 &2 &2 &1 &1 &1 &4 &7 &1 \\
        \hline
        \hline
    \end{tabular}
    \begin{flushleft}
        \footnotesize{$\# 1$ = $|\mathcal{A}_R| - |A_*||\mathcal{K}| - |A_*|(|A_*|-1)$} \\
        \footnotesize{$\# 2$ = $|A_*|(|A_*|-1) - (|A_*|-1)$} \\ 
        \footnotesize{$\# 3$ = $|\mathcal{A}_R| + (|A_D| - |V| + 1)|\mathcal{K}||\mathcal{L}| - (|\delta(S)| - r) - \sum^r_{i=1} (|A(S_i)|-|S_i|+1)$} \\
        \footnotesize{$\# 4$ = $\sum^r_{i=1} (|A(S_i)|-|S_i|+1)$} \\ 
    \end{flushleft}
    \vspace{-0.4cm}
\end{table*}

\begin{table*}[!htp]\centering
\begin{tabular}{lp{5cm}}\toprule
\multicolumn{2}{p{6.5in}}{\textbf{Visual proof of Lemma \ref{lem:conn}:} First the set S is identified such that it satisfies the claims from Lemma \ref{lem:conn}. The tags are highlighted one at a time to simplify the construction of an upper triangular matrix. Observe that all the solutions below satisfy the flow, source, and service constraints. In addition, these satisfy the service constraints of $A_q$ with equality i.e. every solution visits only one required arc from the set $A_q$.} \\ 
\hline
\textbf{Affinely independent solution} &\textbf{Description} \\\midrule
\multicolumn{2}{p{6.5in}}{\textbf{Row-3 Column-2 of Table \ref{tab:lbafsol}}} \\ 
\raisebox{-\totalheight}{\includegraphics[width=0.6\textwidth, height=42mm]{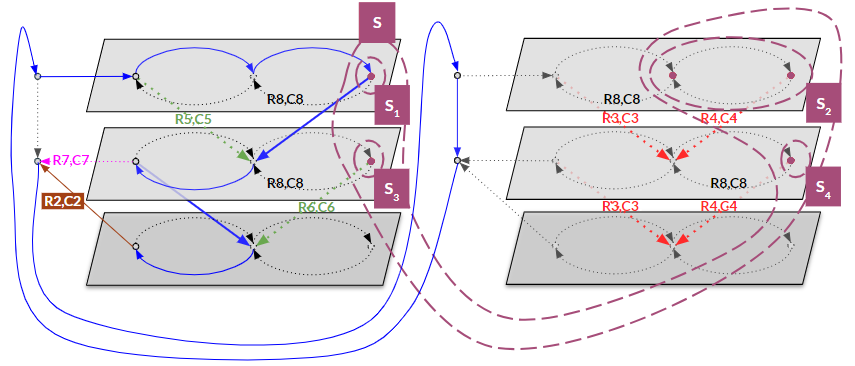}} & 
\begin{itemize}
    \setlength{\parskip}{0pt}
    \setlength{\itemsep}{0pt plus 1pt}
    \item[] Type: Base solution $X_0$ 
    \item[] Highlighted tag: R2,C2 
    \item[] Note:  Set $S = \cup_{r} S_r$ is allowed to be disjoint. Vertices and arcs belonging to the disjoint sets are encircled. 
\end{itemize} \\
& \\ 
\multicolumn{2}{p{6.5in}}{\textbf{Row-4 Column-3 of Table \ref{tab:lbafsol}}} \\ 
\raisebox{-\totalheight}{\includegraphics[width=0.6\textwidth, height=42mm]{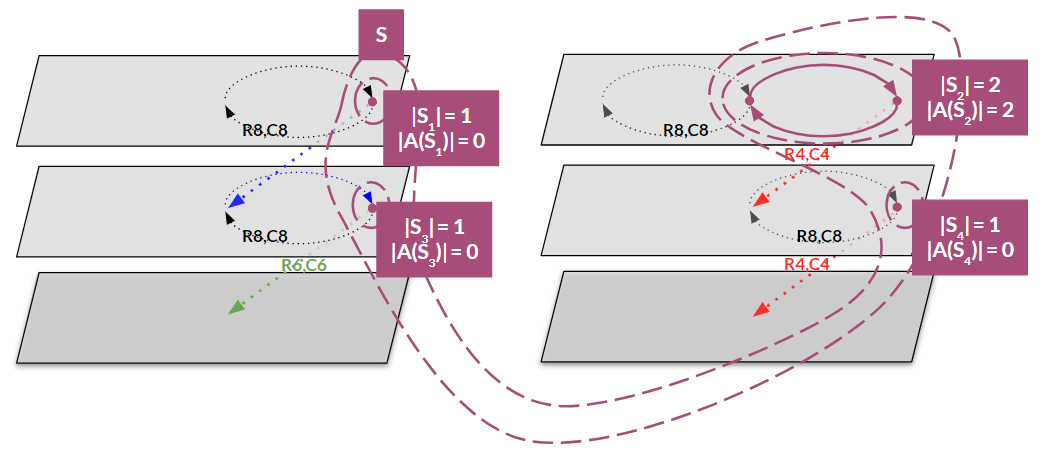}} & 
\begin{itemize}
    \setlength{\parskip}{0pt}
    \setlength{\itemsep}{0pt plus 1pt}
    \item[] Note: The vertex and arc set sizes of the disjoint sets are mentioned in the pictures. The boundary arcs are represented with dotted arrows. 
\end{itemize} \\
\raisebox{-\totalheight}{\includegraphics[width=0.6\textwidth, height=42mm]{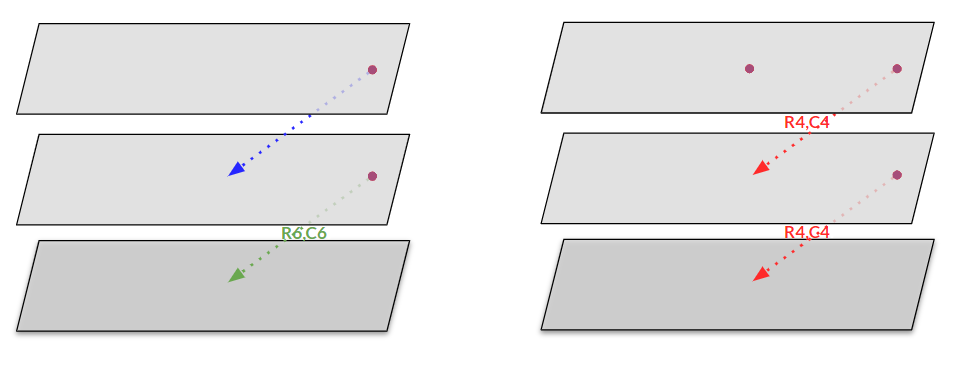}} & 
\begin{itemize}
    \setlength{\parskip}{0pt}
    \setlength{\itemsep}{0pt plus 1pt}
    \item[] Note: The copy of replicated arcs (say $A_q$) should be contained in the arcs set of S or its neighbors i.e. $A_q$ is a subset of $A(S) \cup \delta(S)$. In this illustration, the set $A_q$ is contained in the boundary arcs' set. 
\end{itemize} \\
\raisebox{-\totalheight}{\includegraphics[width=0.6\textwidth, height=42mm]{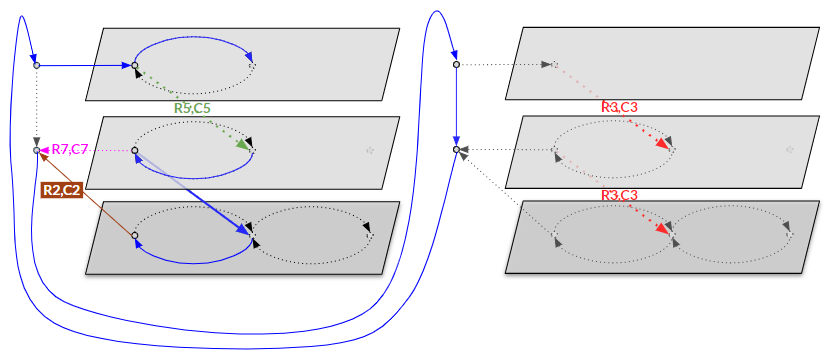}} & 
\begin{itemize}
    \setlength{\parskip}{0pt}
    \setlength{\itemsep}{0pt plus 1pt}
    \item[] Note: Part of the base solution $X_0$ represented over the graph $\mathcal{G}(\mathcal{V}\backslash S)$.  
\end{itemize} \\
\bottomrule
\end{tabular}
\end{table*}

\newpage
\begin{table*}[!htp]\centering
\begin{tabular}{lp{5cm}}\toprule
\textbf{Affinely independent solution} &\textbf{Description} \\\midrule
\multicolumn{2}{p{6.5in}}{\textbf{Row-2 Column-2 of Table \ref{tab:lbafsol}, $\# solutions = 1$ ($1$ of $r$)}} \\ 
\raisebox{-\totalheight}{\includegraphics[width=0.6\textwidth, height=42mm]{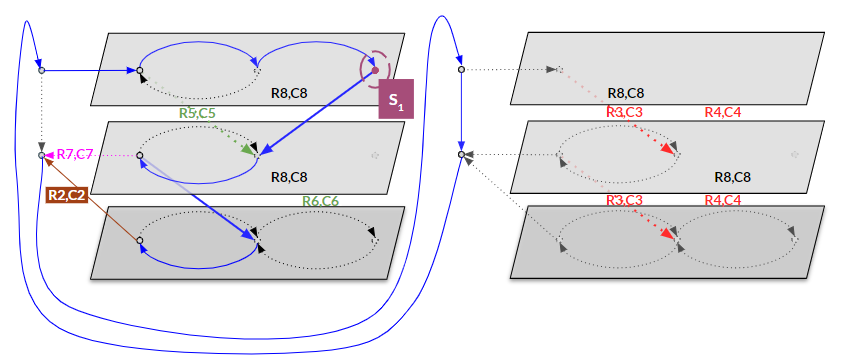}} & 
\begin{itemize}
    \setlength{\parskip}{0pt}
    \setlength{\itemsep}{0pt plus 1pt}
    \item[] Note: The base solution visits set $S_1$. Hence, one of the $r$ solutions is constructed, included in Step 2 of the proof that generates $r$ solutions based on arcs that go in and out.  
\end{itemize} \\
\multicolumn{2}{p{6.5in}}{\textbf{Row-3 Column-3 of Table \ref{tab:lbafsol},  $\# solutions = |\mathcal{A}_R| - |A_*||\mathcal{K}| - |A_*|(|A_*|-1) = 2$}} \\
\multicolumn{2}{p{6.5in}}{\textbf{   (all required arc - required arcs in $A_q$ - other required arcs of $1^\textit{st}$ agent-sub-graph)}} \\ 
\raisebox{-\totalheight}{\includegraphics[width=0.6\textwidth, height=42mm]{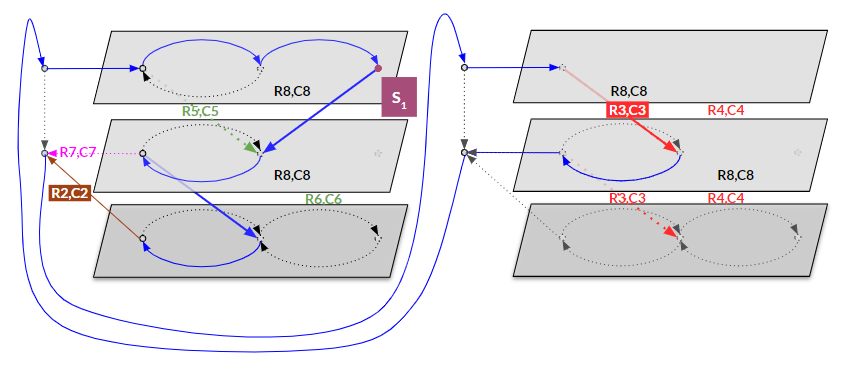}} & 
\begin{itemize}
    \setlength{\parskip}{0pt}
    \setlength{\itemsep}{0pt plus 1pt}
    \item[] Type: Using base solution $X$
    \item[] Highlighted tag: R3,C3 (1)
    \item[] Note: These affinely independent solutions are accounted for in the 1st paragraph of Step-1 of the proof. 
\end{itemize} \\
\raisebox{-\totalheight}{\includegraphics[width=0.6\textwidth, height=42mm]{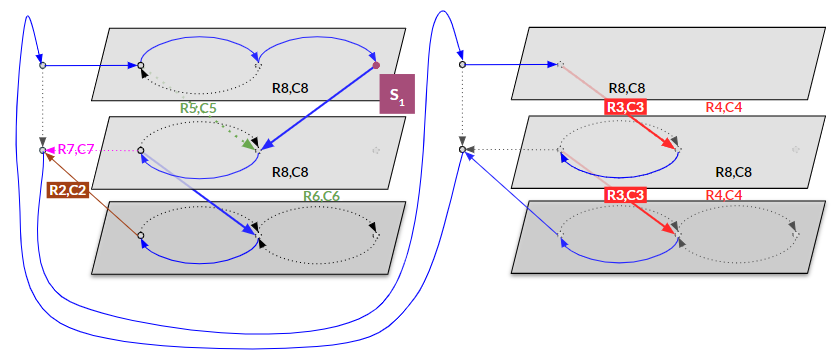}} & 
\begin{itemize}
    \setlength{\parskip}{0pt}
    \setlength{\itemsep}{0pt plus 1pt}
    \item[] Type: Using base solution $X$
    \item[] Highlighted tag: R3,C3 (2)
\end{itemize} \\
\multicolumn{2}{p{6.5in}}{\textbf{Row-4 Column-4 of Table \ref{tab:lbafsol},  $\# solutions = |A_*||\mathcal{K}| - |A_*| = 2$, ($2~\&~3$ of $r$)}} \\
\multicolumn{2}{p{6.5in}}{\textbf{   (required arcs in $A_q$ - required arcs of $1^\textit{st}$ agent-sub-graph of $A_q$)}} \\ 
\raisebox{-\totalheight}{\includegraphics[width=0.6\textwidth, height=42mm]{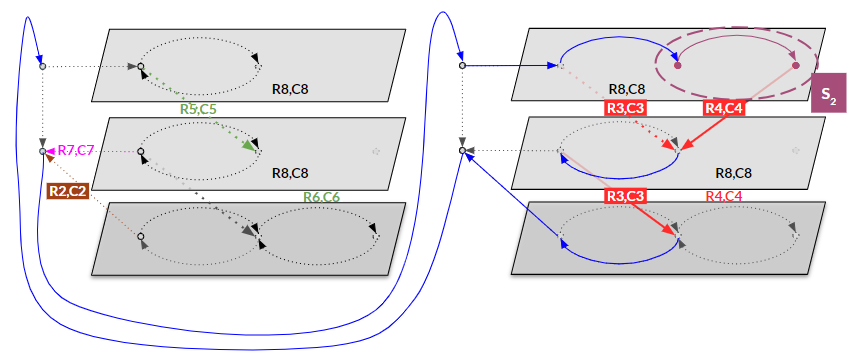}} & 
\begin{itemize}
    \setlength{\parskip}{0pt}
    \setlength{\itemsep}{0pt plus 1pt}
    \item[] Type: Base solution $X_0'$
    \item[] Highlighted tag: R4,C4 (1)
    \item[] Note: These affinely independent solutions are also accounted for in the 1st paragraph of Step-2 of the proof. The required arcs from $A_q$, not present in the $1^\textit{st}$ agent-sub-graph are identified in R4,C4. 
\end{itemize} \\
\bottomrule
\end{tabular}
\end{table*}

\newpage
\begin{table*}[!htp]\centering
\begin{tabular}{lp{5cm}}\toprule
\textbf{Affinely independent solution} &\textbf{Description} \\\midrule
\raisebox{-\totalheight}{\includegraphics[width=0.6\textwidth, height=42mm]{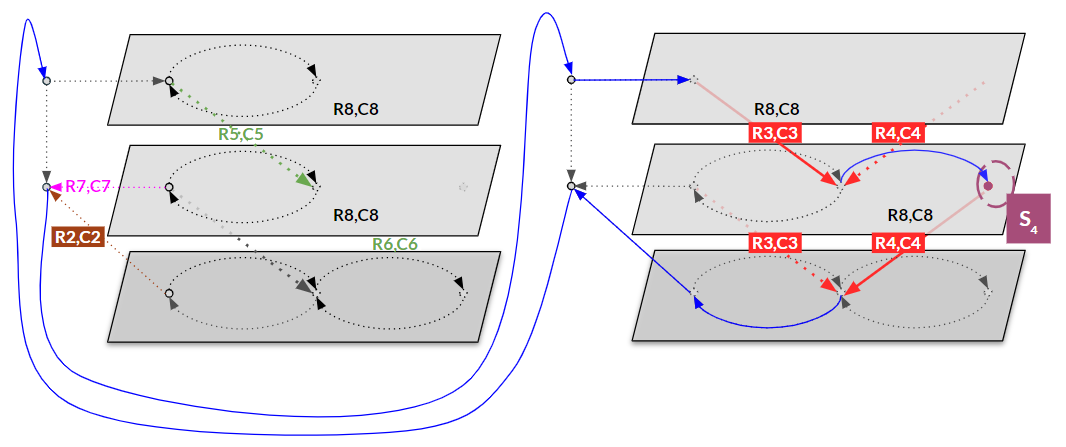}} & 
\begin{itemize}
    \setlength{\parskip}{0pt}
    \setlength{\itemsep}{0pt plus 1pt}
    \item[] Type: Using base solution $X_0'$
    \item[] Highlighted tag: R4,C4 (2)
\end{itemize} \\
\multicolumn{2}{p{6.5in}}{\textbf{Row-5 Column-5 of Table \ref{tab:lbafsol},  $\# solutions = |A_*|(|A_*|-1) - (|A_*|-1) = 1$}} \\
\multicolumn{2}{p{6.5in}}{\textbf{   (required arcs not in $A_q$ of the $1^\textit{st}$ agent-sub-graph - required arcs in $X_0$ among these)}} \\ 
\raisebox{-\totalheight}{\includegraphics[width=0.6\textwidth, height=42mm]{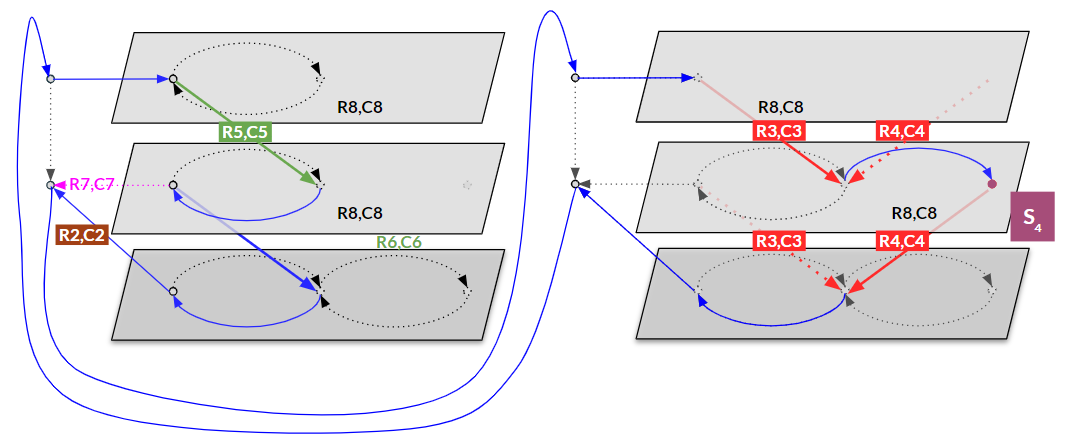}} & 
\begin{itemize}
    \setlength{\parskip}{0pt}
    \setlength{\itemsep}{0pt plus 1pt}
    \item[] Type: Using base solution $X_0'$
    \item[] Highlighted tag: R5,C5
\end{itemize} \\
\multicolumn{2}{p{6.5in}}{\textbf{Row-6 Column-6 of Table \ref{tab:lbafsol},  $\# solutions = (|A_*|-1) = 1$ ($4$ of $r$)}} \\
\multicolumn{2}{p{6.5in}}{\textbf{   (required arcs neither in $A_q$ of the $1^\textit{st}$ agent-sub-graph nor in $X_0$)}} \\ 
\raisebox{-\totalheight}{\includegraphics[width=0.6\textwidth, height=42mm]{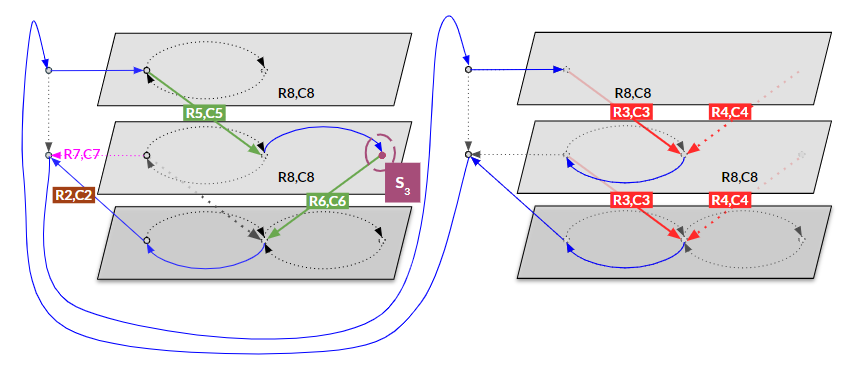}} & 
\begin{itemize}
    \setlength{\parskip}{0pt}
    \setlength{\itemsep}{0pt plus 1pt}
    \item[] Type: Using base solution $X_0'$
    \item[] Highlighted tag: R6,C6 
    \item[] Note: This is another affinely independent solution that visits one of the set $S_r$. This is accounted in the Step-2 first paragraph - one among the $r$ solutions. 
\end{itemize} \\
& \\ 
\multicolumn{2}{p{6.5in}}{\textbf{Row-7 Column-7 of Table \ref{tab:lbafsol},  $\# solutions = |A_*| - 1 = 1$}} \\
\multicolumn{2}{p{6.5in}}{\textbf{   (depot-sink arcs except from the last layer of $1^\textit{st}$ agent-sub-graph)}} \\ 
\raisebox{-\totalheight}{\includegraphics[width=0.6\textwidth, height=42mm]{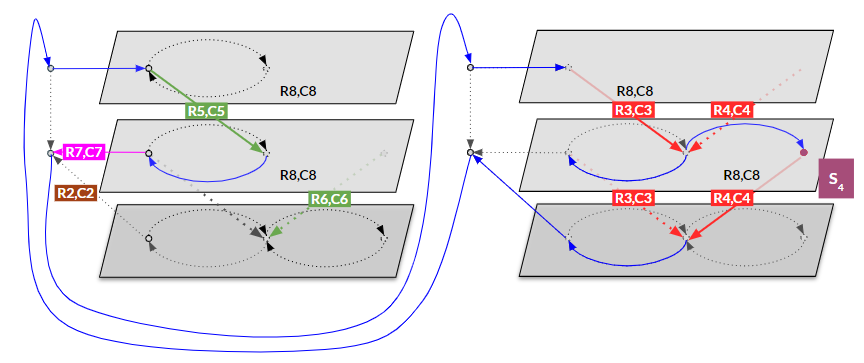}} & 
\begin{itemize}
    \setlength{\parskip}{0pt}
    \setlength{\itemsep}{0pt plus 1pt}
    \item[] Type: Using base solution $X_0'$
    \item[] Highlighted tag: R7,C7 
    \item[] Note: Utilizing the depot-sink arcs of the $1^\textit{st}$ agent-sub-graph, that were not taken in the earlier construction of affinely independent solutions. 
\end{itemize} \\
\bottomrule
\end{tabular}
\end{table*}

\newpage
\begin{table*}[!htp]\centering
\begin{tabular}{lp{5cm}}\toprule
\textbf{Affinely independent solution} &\textbf{Description} \\\midrule
\multicolumn{2}{p{6.5in}}{\textbf{Row-8 Column-8 of Table \ref{tab:lbafsol},  $\# solutions = \delta(S)-2r = 4$}} \\
\multicolumn{2}{p{6.5in}}{\textbf{   (arcs from the boundary - arcs already chosen from the boundary)}} \\ 
\raisebox{-\totalheight}{\includegraphics[width=0.6\textwidth, height=42mm]{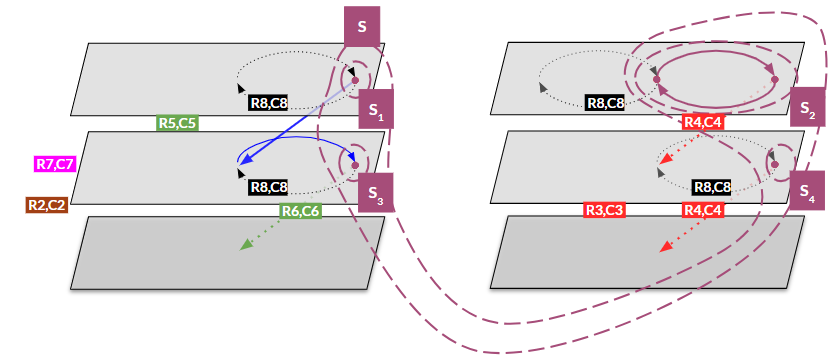}} & 
\begin{itemize}
    \setlength{\parskip}{0pt}
    \setlength{\itemsep}{0pt plus 1pt}
    \item[] Type: Using base solution 
    \item[] Type: Using base solution $X_0$ and $X_0'$
    \item[] Highlighted tag: R8,C8 
    \item[] Note: For visiting each of the $r$ sets of $S_r$, only two of the boundary arcs were utilized, therefore $2r$ is removed from the total set of boundary arcs to construct the affinely independent solutions.   
\end{itemize} \\
\multicolumn{2}{p{6.5in}}{\textbf{Row-9 Column-9 of Table \ref{tab:lbafsol},  $\# \textit{cycles in }\mathcal{G}(\mathcal{V}\backslash S) = (|A_D| - |V| + 1)|\mathcal{K}||\mathcal{L}| - ((|\delta(S)| - r) + \sum^r_{i=1} (|A(S_i)|-|S_i|+1) - |A_*||\mathcal{K}|) = 7$}} \\
\multicolumn{2}{p{6.5in}}{\textbf{   (all intra-layer cycles constructed using the graph $\mathcal{G}(\mathcal{V}\backslash S)$)}} \\ 
\multicolumn{2}{p{6.5in}}{\textbf{   = all intra-layer cycles - (boundary arcs contributing to cycles + cycles in the set $S$)}} \\ 
\raisebox{-\totalheight}{\includegraphics[width=0.6\textwidth, height=42mm]{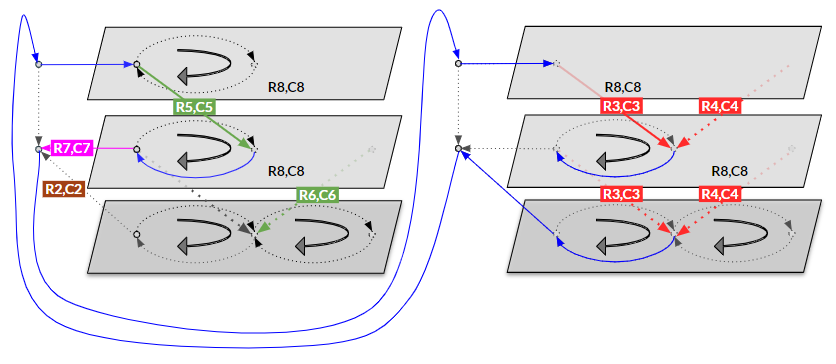}} & 
\begin{itemize}
    \setlength{\parskip}{0pt}
    \setlength{\itemsep}{0pt plus 1pt}
    \item[] Type: Using base solution $X_0$ and $X_0'$
    \item[] Note: The cycles act as extreme rays. Any integer multiple times the extreme ray/cycle can be incorporated into a suitable solution constructed earlier. These cycles based solutions are considered in the $2^\textit{nd}$ paragraph of Step-1 of the proof. 
\end{itemize} \\
\multicolumn{2}{p{6.5in}}{\textbf{Row-10 Column-10 of Table \ref{tab:lbafsol},  $\# \textit{cycles in }\mathcal{G}(S) = \sum^r_{i=1} (|A(S_i)|-|S_i|+1) = 1$}} \\
\multicolumn{2}{p{6.5in}}{\textbf{   (all intra-layer cycles constructed using the graph $\mathcal{G}(\S)$)}} \\ 
\raisebox{-\totalheight}{\includegraphics[width=0.6\textwidth, height=42mm]{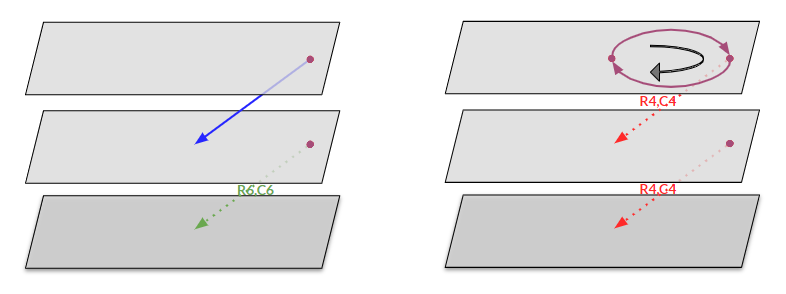}} & 
\begin{itemize}
    \setlength{\parskip}{0pt}
    \setlength{\itemsep}{0pt plus 1pt}
    \item[] Type: Using base solution $X_0$ and $X_0'$
    \item[] Note: These cycles based solutions are considered in the $2^\textit{nd}$ paragraph of Step-2 of the proof. 
\end{itemize} \\
\bottomrule
\end{tabular}
\end{table*}